\documentclass[article,ij4uq]{ij4uq}              %  final version
\RequirePackage{snapshot}
\usepackage{soul}
\usepackage{comment}
\newcommand{\new}[1]{#1}
\newcommand{\old}[1]{#1}
\renewcommand{\old}[1]{}

\usepackage[hang]{footmisc}
\fancypagestyle{plain}{%
  \fancyhf{}
  \fancyhead[R]{\small {\it International Journal for Uncertainty Quantification}, x(x): \thepage--\pageref{LastPage} (\myyear\today)}
  \fancyfoot[R]{\small\bf\thepage }
  \fancyfoot[L]{\fottitle}
  }

\usepackage{amsmath, amsfonts, bm, bbm, graphicx, caption, subcaption, url}
\usepackage[nomarkers,figuresonly]{endfloat} % Use this to force all figures to the end of the document
\newcommand*{\vertbar}{\rule[-1ex]{0.5pt}{2.5ex}}

\begin{document}

\volume{Volume x, Issue x, \myyear\today}
\title{Model Structural Inference using Local Dynamic Operators}
\titlehead{Model Structural Inference using LDOs}
\authorhead{A.M. DeGennaro, N.M. Urban, B.T. Nadiga, \& T. Haut}
%For at least  authors with different addresses, use instead the following commands
\corrauthor[1]{Anthony M. DeGennaro}
\author[2]{Nathan M. Urban}
\author[2]{Balasubramanya T. Nadiga}
\author[3]{Terry Haut}
\corremail{adegennaro@bnl.gov}
\corraddress{Computational Science Initiative, Brookhaven National Laboratory, Upton, NY, 11973}
\address[1]{Computational Science Initiative, Brookhaven National Laboratory, Upton, NY, 11973}
\address[2]{Computer, Computational, and Statistical Sciences, Los Alamos National Laboratory, Los Alamos, NM, 87544}
\address[3]{Computational Physics Group, Lawrence Livermore Laboratory, Livermore, CA, 94550}
% End information for at least  authors with different addresses
% For authors with the same post address,
%\corrauthor{First A. Author}
%\corremail{f.author@affiliation.com}
%\author{Second B. Author, Jr.}
%\address{Department of Chemistry and Courant, Institute of Mathematical Sciences, New York, NY 10012, USA}
% End commands for all authors with the same address

\dataO{mm/dd/yyyy}
%\dataF{mm/dd/yyyy}

\abstract{This paper focuses on the problem of quantifying the effects
  of model-structure uncertainty in the context of time-evolving
  dynamical systems. This is motivated by multi-model uncertainty in
  computer physics simulations: developers often make different
  modeling choices in numerical approximations and process
  simplifications, leading to different numerical codes that
  ostensibly represent the same underlying dynamics. We consider
  model-structure inference as a two-step methodology: the first step
  is to perform system identification on numerical codes for which it
  is possible to observe the full state; the second step is structural
  uncertainty quantification (UQ), in which the goal is to search
  candidate models ``close'' to the numerical code surrogates for
  those that best match a quantity-of-interest (QOI) from some
  empirical dataset. Specifically, we: (1) define a discrete, local
  representation of the structure of a partial differential equation,
  which we refer to as the ``local dynamical operator'' (LDO); (2)
  identify model structure non-intrusively from numerical code output; (3)
  non-intrusively construct a reduced order model (ROM) of the
  numerical model through POD-DEIM-Galerkin projection; (4) perturb
  the ROM dynamics to approximate the behavior of alternate model
  structures; and (5) apply Bayesian inference and energy conservation
  laws to calibrate a LDO to a given QOI. We demonstrate these
  techniques using the two-dimensional rotating shallow water (RSW)
  equations as an example system.}


\maketitle

\section{Introduction}
For predictive modeling of a complex dynamical system to be
successful, it is critical to be able to comprehensively account for
the full range of uncertainties that influence the system. Such
uncertainties are found in initial conditions, boundary and/or
operating conditions, governing parameters, computational approaches,
and model-form. The last of these stems from an imperfect or
incomplete description of the full range of dynamics and/or physics
that underlies the complex behavior exhibited by the system of
interest.

While model-form uncertainty is an active area of research, it is at
present less studied than other forms of uncertainty. A main approach
to quantifying it has relied on multi-model ensembles or ensembles of
opportunity wherein models from different modeling groups are analysed
statistically in an {\it a posteriori} fashion~\cite{solomon2007climate,tebaldi2005quantifying,krishnamurti2000multimodel}. However, this approach
typically does not give representation to all plausible structural
differences: there are many hypothetical model structures,
representing reasonable physical/empirical assumptions that in
principle could be implemented and studied, but have not been. A
multi-model ensemble approach therefore generally cannot quantify the
full range of possible structural uncertainties in a dynamical system,
since it studies only a discrete subset of model structures.

In contrast, parametric uncertainty is a well-studied canonical
problem, with correspondingly well-developed techniques for
quantification. Motivated by this observation, one might wonder if it
is possible to effectively {\it re-cast} certain problems in
structural uncertainty as equivalent problems in parametric
uncertainty, for which a broad and powerful set of tools exist. This
is a main component of the approach we propose in this paper: we seek
to encapsulate plausible differences in model structure in some
parametric approximation of the governing dynamics. In this scheme,
different models are represented by different locations in some
parameter space.

%% In order to do this, it is necessary to first define a ``model
%% structure space'', analogous to model parameter space, that formalizes
%% the modeling choices leading to different structures. The source code
%% of a model itself is a representation of the model structure, but this
%% can run to millions of lines and is not easily manipulable from a
%% practical perspective.  We seek a relatively low-dimensional,
%% non-intrusive means to locate a particular model within ``structure
%% space''.

In order to specify the nature of this structure-parameterization, we
begin by making the crucial assumption that the field dynamics may be
well approximated by a functional relationship which is {\it local}
and {\it spatio-temporally invariant}. By local, we mean that the field
evolution in time at a particular point in space is governed only by
the field values in a nearby neighborhood; by spatio-temporal
invariance, we mean that the mapping that describes this evolution is
the same for all points in space at all times.

Given these assumptions, the discretized dynamics can be described by
what we will refer to as a local dynamic operator (LDO), which is
simply a functional relationship between spatially-local field values
that approximates the discretized governing field dynamics at a
spatial point. For example, if the governing equations are hyperbolic,
then the LDO is a function that takes field values in a
spatially-local neighborhood of a center point and outputs the field
value at that center point, one time step forward in time. \old{We
  should also emphasize that a LDO -- being some discrete
  approximation to the underlying field equations -- is usually not
  uniquely defined for a given set of dynamics.} \new{Note that there
  is an attractive consequence of our assumptions of locality and
  spatio-temporal invariance with respect to system identification: if
  we wish to infer a LDO from numerical/experimental data, access to
  the full global state vector is not strictly required. We may simply
  collect data from a subset of spatial points (together with the
  appropriate surrounding local neighborhoods). This is a notable
  advantage relative to a system identification technique that would
  require the full global state vector (e.g., POD).}

As all of the dynamics are encoded in the LDO, any structural
uncertainties are as well. Furthermore, if we can design the LDO to be
a weighted sum of different elementary functions of the local field
values, then the relevant structural uncertainties manifest themselves
as uncertainties in the values of the weights, which are simply
parameters. This is a sketch of the process by which we convert
structural uncertainties to parametric ones. \old{Note that the LDO,
  being a local operator, is independent of the PDE domain size or
  mesh resolution. This inherent low-dimensionality is advantageous if
  we wish to infer a LDO from numerical/experimental data.}

Having formalized a means to parameterize model structure, all of the
machinery of parametric uncertainty quantification (UQ) is available
to study the structural uncertainty in model dynamics, which we
propose doing in the following \old{way} \new{two-step process}. First, we show how we
may learn model structure (i.e., LDO parameters) non-intrusively
using statistical function approximation techniques (e.g.,
least-squares regression). Second, we show how Monte Carlo sampling
may be used to perturb the LDO parameters, giving us a systematic,
probabilistic means for performing model structural UQ. \old{Next, we
demonstrate how physical constraints (e.g., conservation laws) can be
used to guide the sampling process by ruling out implausible,
non-physical model structures. Finally, we investigate building
reduced-order models (ROMs) of the LDO-parameterized dynamics, using
POD-DEIM-Galerkin projection. The motivation for this is efficiency: a
ROM simulation is much less computationally expensive than a
simulation of the full dynamics. We examine whether a ROM identified
for one model structure (i.e., set of LDO parameters) is of sufficient
generality to reasonably approximate the dynamics of perturbed models.}

\new{Without any other constraints, the number of LDO parameters can
  be prohibitively large for the second step of our process: Monte
  Carlo sampling in a high-dimensional parameter space is a
  challenge. We explore two methods for alleviating this. First, we
  explicitly assume that the user has partial prior knowledge about
  the physics involved, and we use this knowledge to constrain the LDO
  parameter space and reduce its dimensionality. In doing this, we are
  making our mehod a ``greybox'', and as such, our approach is not
  aimed at users who desire a purely data-driven, completely
  physics-agnostic approach. Instead, we are targeting those users who
  have a combination of both some prior knowledge and some data, and
  we are proposing a framework to effectively join the two. The second
  avenue that we investigate for reducing computational expense
  involves building reduced-order models (ROMs) of the
  LDO-parameterized dynamics, using POD-DEIM-Galerkin projection.}

The envisioned practical application for these ideas is
quantification/learning of subgrid physics. The motivation here is
that several numerical models exist that each close the subgrid
physics for a particular physics problem (e.g., climate simulations),
but each of these individual models have been developed using
different physical assumptions and/or empirical data and so produce
different long-time integrated behaviors. However, if we are able to
approximate the input-output behavior of these different models using
a dictionary of local state features, then we can learn that area of
the LDO feature space that contains those models and search in some
local neighborhood of that for candidate models that best reproduce a
quantity-of-interest (QOI) taken from some dataset that we regard as
the ``truth''. This can be thought of as a two-step process of system
identification followed by structural UQ: the first step is to learn
feature space representations of the dynamics of several numerical
models (system identification); the second step is to search perturbed
candidate models ``near'' these numerical models for those that best
match empirical data QOIs (structural UQ).

This paper builds on prior work in non-intrusive system identification
and also extends those ideas to the setting of structural UQ as well. We
briefly review some of the popular existing techniques for comparison
here. Proper Orthogonal Decomposition (POD) and POD-Galerkin
projection~\cite{Lumley,HolmesBook,SirovichPOD,BerkoozPOD} is perhaps
the most widely used \old{non-}intrusive model reduction technique; its goal
is to construct a set of basis vectors for the state from data, and
then to project the governing equations onto that empirical basis in a
Galerkin fashion. Our approach, in contrast, uses local functions of
the state rather than the global state, and attempts to learn
equations for the dynamics, as opposed to just global spatial
correlations; also, we need not know the underlying dynamical
equations, as is required for POD-Galerkin projection.

Dynamic Mode Decomposition~\cite{RowleyDMD,SchmidDMD,TuDMD,
  WilliamsEDMD, WilliamsKDMD, Hemati2014, ProctorDMD} is another
popular model inference technique; it approximates the Koopman
operator for a dynamical system using data, which can be used to
propagate global functions of the state. This, again, is an approach
that uses the global state vector and does not assume spatial locality
of the dynamics as we do; also, DMD is not an equation learning
technique.

\old{Loiseau and Brunton}\new{Brunton, Proctor, Kutz and
  collaborators} recently introduced another set of non-intrusive
techniques aimed at learning governing dynamics equations, called
constrained sparse Galerkin regression (CSGR) and sparse
identification of nonlinear dynamics
(SINDy)~\cite{Loiseau,SINDY,rudy}. \new{The first step of our
  methodology (i.e., system identification) is essentially equivalent
  to their approach: their method learns governing equations from data
  in some lifted feature space consisting of functions of the state,
  and physical constraints may be enforced as well. However, the final
  goal of our methodology is our second step of structural UQ, in
  which we explore and quantitatively evaluate a large space of novel
  model structures with respect to some empirical QOI, and each of
  these candidate models are perturbed from some ``base''
  model. System identification is just the first step we use to
  identify that base model (as well as the model structural space
  itself) and to inform our prior distribution about that model for
  the structural UQ step.} \old{Our methodology differs from these in
  important ways; most notably, CSGR/SINDy -- like POD-Galerkin
  projection and DMD -- does not assume fundamental locality of PDE
  model dynamics as we do. In assuming that, we introduce an
  attractive reduction of dimensionality of the dynamics that is not
  built in to CSGR/SINDy.  Additionally, CSGR/SINDy is focused more on
  system identification than structural UQ, and the dimension
  independence of our LDO approach may become particularly important
  if the goal is to sample a large space of novel, perturbed model
  structures.}

Peherstorfer~\cite{Peherstorfer,Peherstorfer2} considered the problem
of inferring approximate dynamics equations non-intrusively from data
by regressing quadratic polynomial combinations of the global state
vector against model output. This is, again, a technique that uses the
global state vector, and focuses only on quadratic polynomial
functions of that state. Our technique uses more general functions of
the local state for approximating the dynamics.

Quade~\cite{SymbolicRegression} demonstrated how a surrogate set of
equations for a dynamical system may be learned using machine learning
techniques from a dictionary of elemental state functions. In contrast
to our work, they considered dynamical systems with a small number of
state variables, and therefore did not consider spatio-temporal field
equations -- and the concept of spatial locality of the dynamics --
that form the focus of our paper.

Recently, Sirignano~\cite{Sirignano} applied machine learning
techniques to train a neural network to satisfy a known differential
operator with known initial/boundary conditions, by batch-feeding
randomly sampled data to the network and updating it with stochastic
gradient descent. Their goal, therefore, was to train an emulator that
could correctly predict solutions to a given set of governing
dynamics, not to learn those governing dynamics themselves (which is
one of our goals). Also, their approach did not make use of the
concept of spatial locality, or physical constraints.

\old{The model-inference technique that is perhaps most} \new{Another
  model-inference technique that is} closely related to ours is the
``blackbox stencil interpolation method'' (BSIM) introduced by
Chen~\cite{ChenThesis}. This method attempts model inference by
regressing model output at a spatial location against the field values
on the stencil at that location. BSIM focuses on approximating the
discretized dynamics with neural network function regression, and
building ROMs for the resulting models. We seek to extend these ideas
in a number of \old{novel} ways. First, we introduce specific
assumptions about the form of the discretized field dynamics by
parameterizing it in terms of a dictionary of functions of the local
state variables. This is in contrast to a more general machine
learning approach in which one need not assume anything about the
functional form of the dynamics. Introducing this parameterization
will allow us to assume a basic structure to the LDO, which makes
trivial the possibility of further constraining this structure with
physical constraints (e.g., energy conservation laws). This
\old{contribution} serves as a means by which the dimensionality of
the LDO parameter space may be reduced by only considering those LDOs
that respect certain physics. \new{Most importantly, BSIM is purely a
  system identification technique, and as such, we reiterate that our
  contribution is a two-step process of system identification used in
  service of structural UQ.} \old{We further note that the main focus
  of BSIM is system identification, not structural uncertainty
  quantification, which is one of our main focuses.}

%% CSGR also performs system identification and constructs ROMs from the
%% identified dynamics.  Unlike BSIM but like our LDO approach, it
%% exploits physical constraints on the identified dynamics, such as
%% energy conservation.  However, unlike the BSIM and LDO approaches, it
%% does not assume the fundamental \textit{locality} of PDE model
%% dynamics.  Its system identification method, sparse identification of
%% nonlinear dynamics (SINDy), performs regression on the global state
%% vector defined over the entire computational domain, using dimension
%% reduction (basis projection and truncation) to deal with scalability
%% issues.  CSGR is thus more general, and while it exploits sparsity in
%% the dynamics, it does not directly exploit the essential dimension
%% independence of the dynamics implied by this class of local problems.
%% Like BSIM, CSGR is also focused more on system identification than
%% structural uncertainties, and the dimension independence of the LDO
%% (or BSIM) approach may become particularly important if the goal is
%% to sample a large space of novel, perturbed model structures.

This paper is organized as follows. In the first section, we begin
with basic definitions of a LDO and related concepts, and then
restrict ourselves to studying only LDOs that have a certain
structure. Next, we cover several strategies for inferring a LDO,
given different assumptions about the simulation output that we are
able to observe. As the space of LDO parameters can be prohibitively
large for inference, we then introduce physical constraints and show
how they may be used to reduce the dimensionality of the allowable LDO
parameter space. In order to further reduce the computational burden
associated with LDO inference, we next propose to use a reduced-order
model (ROM) as a surrogate for the LDO in inference problems. Finally,
we conclude with several numerical examples, using the rotating
shallow water (RSW) equations as the base dynamics. We conclude, based
on these examples, that LDO inference using observable simulation
output is possible, and can be accelerated with some success by using
a ROM for the LDO.

\section{Local Dynamic Operators}\label{sec:LDO}
We begin the discussion by making precise our definitions of a Local
Dynamic Operator (LDO). This is the means by which we will capture
model-form uncertainties in a parametric setting. We define the LDO as
a discrete approximation, in some local feature space, to some
relevant continuum governing dynamics. As a general exposition,
consider the field equations:
\begin{equation}
\bm{\dot{u}} = \bm{F}(\bm{u}) \; ,
\label{eq:PDE}
\end{equation}
where $\bm{u}(\bm{x},t) \in \mathbb{R}^N$ is the state \new{(i.e.,
  there are $N$ state variables defined at a given point in space and
  time)}, $\bm{F}(\bm{u}) : \mathbb{R}^N \mapsto \mathbb{R}^N$ is the
continuous-time evolution operator, $\bm{x}$ denotes space, and $t$
denotes time. \new{For clarity, we remark that $\bm{F}(\bm{u})$
  denotes the field equations for the state at a given point in
  space/time, as opposed to the governing equations for the full
  dynamical system.} We assume that there exists a discrete,
spatially-local evolution operator that approximates the governing
equations:
\begin{equation}
\bm{\frac{\Delta u}{\Delta t}} = \bm{f}(\mathcal{L}(\bm{u},\bm{x})) \; .
\label{eq:LDO}
\end{equation}
Here, $\bm{f(\cdot)}$ denotes a discretization of $\bm{F(\cdot)}$, and 
$\mathcal{L}(\cdot)$ is general notation for an operator that produces
the field values in a neighborhood local to a particular point in
space/time. %% For example, this could restrict the spatial domain to
%% only those points within a ball of radius $R$:
%% \begin{equation}
%%     \mathcal{L}(\bm{u},\bm{x},t) =  
%% \begin{cases}
%%     \bm{u}(\bm{x},t) , & \text{if } \bm{|x|}\leq \textit{R}\\
%%     0,                 & \text{otherwise}
%% \end{cases}
%% \end{equation}
For example, in a structured-grid setting, this definition of
$\mathcal{L}(\cdot)$ could effectively produce the field values on the stencil for a
particular point in space/time, i.e.
\begin{equation}
    \mathcal{L}(\bm{u},\bm{x}) = \bm{S} = \begin{bmatrix}
           u_1^{(1)}          \\
           \vdots            \\
           u_N^{(M)}           \\
    \end{bmatrix} \; ,
\end{equation}
where there are $M$ stencil points and $N$ state variables, and
$u_i^{(j)}$ denotes the value of the state variable $i$ at the
$j^{th}$ stencil point. It should be clear that $\bm{f}(\cdot)$ in
Eq.~\ref{eq:LDO} is an operator that acts on a neighborhood of points
(i.e., the stencil $\bm{S}$) local to a specific point. We assume
going forward that $\mathcal{L}(\bm{u},\bm{x}) =
\bm{S}$. \new{Please also note that for the sake of notational
  coherence, both $\bm{u}$ and $\bm{f}$ will be treated as row
  vectors in what follows (i.e., $\bm{u},\bm{f} \in \mathbb{R}^{1 \times N}$).}

Beyond spatial locality, we wish to make further assumptions on the
structure of $\bm{f}(\cdot)$ that will make it tractable to
investigate parametric perturbations. First, let us assume that
$\bm{f}(\cdot)$ may be written as a linear combination of $Q$ basis
functions:
\begin{equation}
\begin{aligned}
\bm{f}(\bm{S}) &= \sum_{i=1}^{Q} \bm{c_i} \psi_i(\bm{S}) \; .
\end{aligned}
\label{eq:LDObasissummation}
\end{equation}
Here, $\psi_j(\cdot) : \mathbb{R}^{NM} \mapsto \mathbb{R}$ is a basis
function with associated coefficient $\bm{c_j} \in
\mathbb{R}^{1 \times N}$. Next, we will incorporate {\it a priori} knowledge of
the types of nonlinearities present in Navier-Stokes to guide our
selection of features, since the examples we consider in this paper
are drawn from fluid dynamics. Thus, in this work, we consider two
specific choices for the basis functions in
Eq.~\ref{eq:LDObasissummation}. One choice is to assume that
$\bm{f}(\cdot)$ may be written as a linear combination of all
polynomial combinations of the local field values up to second order,
which we may express as:
\begin{equation}
\begin{aligned}
\bm{f}(\bm{S}) &= \bm{p_0} + \sum_{i=1}^{NM} \bm{p_i} S_i + \sum_{j=1}^{NM} \sum_{k=j}^{NM} \bm{p_{j,k}} S_jS_k \\
             &= \bm{\psi} \bm{P} \; .
\end{aligned}
\label{eq:LDOsummation}
\end{equation}

%% where:

%% \begin{equation}
%%     \bm{\Psi}(\bm{S}) = \begin{bmatrix}
%%            1                \\
%%            u_1^{(1)}          \\
%%            \vdots            \\
%%            u_N^{(M)}           \\
%%            u_1^{(1)}u_1^{(1)}    \\
%%            \vdots            \\ 
%%            u_N^{(M)}u_N^{(M)}
%%          \end{bmatrix}
%% \;\; , \;\;
%%     \bm{P} = \begin{bmatrix}
%%            \horzbar     &  \bm{p_0}   &   \horzbar              \\
%%            &  \vdots   &                  \\
%%            \horzbar     &  \bm{p_{NM}}   &   \horzbar              \\
%%            \horzbar     &  \bm{p_{1,1}}   &   \horzbar        \\
%%            &  \vdots   &                  \\
%%            \horzbar     &  \bm{p_{NM,NM}}   &   \horzbar               \\
%%          \end{bmatrix}
%% \label{eq:LDOElements}
%% \end{equation}
Here, $\bm{P} \in \mathbb{R}^{Q \times N}$ is the matrix of
coefficients\new{,} \old{and} $\bm{\psi} \in \mathbb{R}^{1 \times Q}$ is the vector of
basis functions, with $Q = 1 + NM + \frac{1}{2}(NM)(NM+1)$\old{.}\new{, $\bm{p_a \in \mathbb{R}^{1 \times N}}$ for any index $a$, and $S_i$ is some element of $\bm{S}$.} Again,
this assumption on the {\it structure} of $\bm{f}(\cdot)$ is chosen to
reflect the structure apparent in many convection-diffusion fluid
equations (e.g., the quadratic nonlinearities in Navier-Stokes).

Of course, Eq.~\ref{eq:LDOsummation} is not the unique way in which
quadratic nonlinearities may be represented. Another (related) choice
is some subset of the basis given above that corresponds to discrete
approximations of select differential operators. One such example
would be:
\begin{equation}
\begin{aligned}
\bm{f}(\bm{u}) &= \bm{p_0} + \sum_{i=1}^{N} \bm{p_{1,i}} u_i^C 
               + \sum_{j=1}^{NX} \bm{p_{2,j}} D_j 
               + \sum_{k=1}^{NX} \bm{p_{3,k}} D^2_k  \\
               + &\sum_{l=1}^{N} \bm{p_{4,l}} (u_l^C)^2 
               + \sum_{m=1}^{N^2X} \bm{p_{5,m}} [\bm{u^C} \otimes D]_m    
               + \sum_{n=1}^{N^2X} \bm{p_{6,n}} [\bm{u^C} \otimes D^2]_n    \\
&= \bm{\psi} \bm{P} \; .
\end{aligned}
\label{eq:LDOsummation2}
\end{equation}
Here, $X$ is the spatial dimension and $D, D^2$ are numerical
approximations to the gradient and Laplacian using the stencil
elements:
\begin{equation}
    D(\bm{S}) = \begin{bmatrix}
           \nabla_{1} (u_1)               \\
           \vdots            \\
           \nabla_{X} (u_N)
         \end{bmatrix}
 \;\;\; , \;\;\;
    D^2(\bm{S}) = \begin{bmatrix}
           \nabla^2_{1} (u_1)               \\
           \vdots            \\
           \nabla^2_{X} (u_N)
         \end{bmatrix} \; .
\end{equation}
\new{As before, $\bm{p_a \in \mathbb{R}^{1 \times N}}$ for any index
  $a$.} Note that $u_j^C$ denotes the center stencil point for state
element $j$ above. This choice of basis is again motivated by the
structure apparent in the quadratic nonlinearities of Navier-Stokes,
but we have additionally limited ourselves to only those terms that
approximate certain differential operators.

Whatever the choice of basis, we refer to the parameters present in
the matrix $\bm{P} \in \mathbb{R}^{Q \times N}$ as the LDO
coefficients. We denote the vector space of all possible LDO
coefficients (equivalently, all coefficient entries possible in the
matrices $\bm{P}$) as $\mathcal{P}$ -- a space of dimension $NQ$. In
all that remains, we will make one of the two previous assumptions on
the \textit{structure} of $\bm{f}(\cdot)$: we will assume it can be
written as in Eq.~\ref{eq:LDOsummation} or as in
Eq.~\ref{eq:LDOsummation2}.

\section{LDO Coefficient Inference}
We now turn our attention to the task of inferring the numerical
values of the LDO coefficients in Eq.~\ref{eq:LDOsummation}. The
intent of this is blackbox inference, wherein one does not know either
the governing dynamics (Eq.~\ref{eq:PDE}) or the numerical
discretization of those dynamics implemented in some complex code, but
nonetheless does have access to input-output data from such a
numerical code, and wishes to learn a LDO that approximates those
dynamics. The LDO coefficients may be estimated numerically from data
through regression, or inferred indirectly through a set of
observables. Which of these techniques is used depends on exactly what
information one has about the dynamics. If one has the ability to
generate and observe snapshots of the full state $\bm{u}(\bm{x},t)$,
then numerical regression of the LDO is possible. However, if
knowledge of neither $\bm{F}(\bm{u})$ nor $\bm{u}(\bm{x},t)$ is
possible, then the LDO coefficients may be inferred indirectly using a
set of observable functions of the state (e.g., via Markov-Chain Monte
Carlo (MCMC) methods~\cite{Metropolis,Hastings}).

\subsection{Data Regression}
\label{sec:dataregress}
If one one has access to high-quality snapshots of the entire state
from numerical simulations, then it is possible to directly infer the
LDO coefficients by simply aggregating training data from different
spatial/temporal points, and then performing some sort of
regression. Assume we collect $N_S$ data snapshots from simulations at
different spatial/temporal locations, and possibly from different
initial conditions. Then, we arrange this data in matrix form:
\begin{equation}
\bm{\dot{u}} \approx \frac{1}{\Delta T} 
    \begin{bmatrix}
           \vertbar & & \vertbar \\
           \bm{\Delta u_1} & \hdots & \bm{\Delta u_N} \\
           \vertbar & & \vertbar 
         \end{bmatrix}
=
\bm{\Psi}\bm{P} \; ,
\label{eq:LDOregression}
\end{equation}
where the matrix columns of $\bm{\Psi} \in \mathbb{R}^{N_S \times Q}$
represent the relevant basis $\bm{\psi}$ evaluated at each of the $N_S$
data snapshots, and $\Delta T$ is the (fixed) timestep between data
snapshots. In Eq.~\ref{eq:LDOregression}, we have approximated
$\bm{\dot{u}}$ with a first-order accurate, forward difference
scheme. Regressing this data (e.g., least-squares, LASSO) yields an
approximation to the matrix of LDO coefficients $\bm{P}$.

\subsection{Bayesian Inference}
If we do not know the PDE describing the system dynamics, and also do
not have the ability to observe the entire state, we may still
estimate the LDO through inference based on observable quantities of
interest (QOIs), which are arbitrary functions of the state (e.g., a
spatio-temporally integrated metric, coarse state data, etc.). This
may be accomplished using a number of methods; in this work, we focus
on the use of the Metropolis algorithm.

We begin by assuming some prior distribution over LDO coefficient
space $\mathcal{P}$, expressing our initial belief about the relative
likelihoods of all possible LDO coefficient combinations. Bayes' law
allows us to revise our statistical belief about the LDO parameters $\bm{P}$ based on measurement of a QOI $q = q(\bm{u},\bm{x},t)$:
\begin{equation}
\pi(\bm{P} | q) = \frac{\pi(q | \bm{P}) \pi_0(\bm{P})}{\pi(q)} \; .
\end{equation}
MCMC sampling can be used to sample this Bayesian posterior distribution (e.g., via the Metropolis algorithm).

\section{Energy Conservation Constraints}\label{sec:Constraints}
One of the goals of the LDO parameterizations we introduced in
\S[\ref{sec:LDO}] is to provide a basis that is both robust enough to
capture a qualitatively ``broad'' range of possible field evolution
dynamics and simultaneously explicit enough that it becomes possible
to constrain the LDO to obey certain physical laws. One reason we may
wish to introduce physical constraints is that we do not wish to
consider field evolution laws that are, in some sense,
non-physical. Another is that introducing additional constraints is
yet another LDO dimension-reduction that could prove helpful,
especially if we wish to infer the LDO coefficients in a Bayesian
setting. Thus, we limit our attention to only those LDOs that satisfy
a statement of energy conservation. As we will show, this process
results in a subspace $\mathcal{P}^E \subset \mathcal{P}$ consisting
of those LDOs that satisfy energy conservation.\old{(see
Fig. 1).} Thus, the goal henceforth is to start
with a statement of energy conservation and find the subspace
$\mathcal{P}^E$ of LDO coefficients that satisfy that
constraint. Structural UQ of the LDO for the new system could then take
place in this constrained space $\mathcal{P}^E$.

\new{The motivating assumption here is that the user has prior
  information regarding some -- but not all -- of the physics
  involved. Intuitively, the user does not know an exact description
  of the physics, but does know that that description -- whatever it
  may be -- must be conformable with some predetermined law. This
  prior knowledge could come from first-principles arguments about the
  problem, e.g. by invoking simplifying assumptions that are
  apparent from the operating conditions or boundary conditions of the
  problem. Regardless, this step will necessarily always be a
  domain-specific process, and it is therefore not possible to offer a
  general calculus for deriving this knowledge. Importantly, it should
  be noted that this partial enforcement of the physics does {\it not}
  imply that some aspects of the physics are more important than
  others: it simply implies that the user {\it knows} those selected
  aspects better than others. The goal of the structural UQ step of
  our algorithm is to learn those residual closures that are not in
  the span of the user's prior knowledge.}

%% \begin{figure}%[h]
%% \centering
%% \includegraphics[width=0.5\textwidth,trim={0 0.0cm 0 0},clip]{EnergySubspace.png}
%% \caption{Generic depiction of $\mathcal{P}$ (the space of all LDO coefficients) and $\mathcal{P}^E$ (the subspace of energy conserving LDO coefficients).}
%% \label{fig:EnergySubspace}
%% \end{figure}

%\subsection{Example 1: Rotating Shallow Water Equations} \label{sec:RSWConstraints}
The numerical examples to come in this paper use the rotating shallow
water equations, and so we demonstrate incorporation of energy
conservation laws using those obeyed by the rapidly rotating shallow
water equations~\cite{Majda1998,MajdaBook}. Note, however, that the
techniques demonstrated here could apply to any differential
conservation statement.

The model rapidly rotating shallow water equations \new{specify the evolution of $N=3$ state variables and} may be written thus:
\begin{equation}
\begin{aligned}
\dot{\bm{u}} + (\bm{u} \cdot \nabla)(\bm{u}) &= 
  \frac{1}{\epsilon} (\hat{\bm{k}} \times \bm{u})
  - \left( \frac{F^{-1/2}}{\epsilon} \right) \nabla \eta \\
\dot{\eta} + \nabla \cdot (\eta \bm{u}) &= 
   - \left( \frac{F^{-1/2}}{\epsilon} \right) \nabla \cdot \bm{u} \; ,
\label{eq:RSW}
\end{aligned}
\end{equation}
where $\bm{u} = (u,v)$ are the $(x,y)$ components of velocity, $\eta$
is the water height, $\hat{\bm{k}}$ is the unit vector normal to the
$(x,y)$ plane, and $F$, $\epsilon$ are model parameters. It may be
shown that this system of equations obeys the following differential
energy conservation law~\cite{Majda1998,MajdaBook}:
\begin{equation}
\frac{\partial{(\eta (E_K + E_P))}}{\partial{t}} 
 + \nabla \cdot (\bm{u} \eta (E_K + 2E_P))
 + (E_K + 2E_P) \nabla \cdot \left(\frac{F^{-1/2}}{\epsilon} \bm{u} \right)
 = 0 \; ,
\label{eq:RSWEnergy}
\end{equation}
where kinetic and potential energies are defined as:
\begin{equation}
\begin{aligned}
E_K = \frac{1}{2}(u^2 + v^2)  \\
E_P = \frac{1}{2} \frac{F^{-1/2}}{\epsilon} \eta \; .
\end{aligned}
\label{eq:EnergyDefn}
\end{equation}

In order to discuss energy constraints for a LDO, we obviously must
first define our LDO basis functions. Let us define them as in
Eq.~\ref{eq:LDOsummation}. We are interested in finding the subspace
of LDO coefficient combinations that satisfy Eq.~\ref{eq:RSWEnergy},
which we will denote as $\mathcal{P}^E \subset \mathcal{P}$. We choose
to use a five-point stencil ($M=5$), and the number of state variables
is $N=3$, so the dimensionality of the unconstrained LDO perturbation
space $\mathcal{P}$ is $3Q = 408$.

One particular solution belonging to $\mathcal{P}^E$ is obviously
given by the rotating shallow water equations themselves, as displayed
in Eq.~\ref{eq:RSW}. Let us denote the LDO corresponding to a
discretization of the shallow water equations as $\bm{P}^{RSW}$. Other
LDOs that conserve energy may be found by examining the null space of
the operator given by the partial time derivative in
Eq.~\ref{eq:RSWEnergy}, which we will denote as $\mathcal{P}^0 \subset
\mathcal{P}$. That is, given LDO coefficients $\bm{P}^{0} \in
\mathcal{P}^0$ that solve the equation:

\begin{equation}
\frac{\partial{(\eta (E_K + E_P))}}{\partial{t}} = 0 \; ,
\label{eq:RSWEnergyNull}
\end{equation}
the linear combination $\bm{P}^E = \bm{P}^{RSW} + \bm{P}^{0}$
will conserve energy.

Expanding Eq.~\ref{eq:RSWEnergyNull} gives:
\begin{equation}
\begin{aligned}
&\dot{\eta}E_K + \eta \dot{E}_K + \dot{\eta} E_P + \eta \dot{E}_P = 0 \\
&\dot{\eta} \left(\frac{1}{2}(u^2 + v^2) \right) + \eta \left( u\dot{u} + v\dot{v} \right)
 + \frac{F^{-1/2}}{\epsilon}\eta \dot{\eta} = 0 \\
&\left(\frac{1}{2}(u^2 + v^2) + \frac{F^{-1/2}}{\epsilon}\eta \right) \bm{\psi} \bm{P}_{\eta}^0
 + (\eta u) \bm{\psi} \bm{P}_{u}^0 + (\eta v) \bm{\psi} \bm{P}_{v}^0  = 0 \; , 
\label{eq:RSWEnergyConstraint}
\end{aligned}
\end{equation}
where the subscripts of $\bm{P}^0$ denote the columns of the
respective state variable. Eq.~\ref{eq:RSWEnergyConstraint} gives
an algebraic set of constraints that must be satisfied by the LDO
parameters $\bm{P}^0$. Examining the structure of these constraints in
concert with the basis given by $\bm{\psi}$ (in Eq.~\ref{eq:LDOsummation})
reveals that the solutions to Eq.~\ref{eq:RSWEnergyConstraint} belong
to the 18-dimensional vector space given by:
\begin{equation}
\begin{aligned}
\bm{\psi} \bm{P}^0 = 
  &\lambda_1 \begin{bmatrix}
           -\frac{1}{2}u^2 - \frac{F^{-1/2}}{\epsilon}\eta        \\
           -\frac{1}{2}uv           \\
           \eta u                      
  \end{bmatrix}^T
+ \lambda_2 \begin{bmatrix}
           -\frac{1}{2}uv           \\
           -\frac{1}{2}v^2 - \frac{F^{-1/2}}{\epsilon}\eta           \\
           \eta v                      
  \end{bmatrix}^T
+ \sum_{i=0}^{15} \lambda_{i+3} \begin{bmatrix}
           -v        \\
           u           \\
           0                      
  \end{bmatrix}^T S_i \; ,
\end{aligned}
\label{eq:RSWNullSpace}
\end{equation}
where $S_0  = 1$ and  $S_1 \ldots S_{15}$  are simply the  elements of
$\bm{S}$  of order  1.
%% + &\lambda_3 \begin{bmatrix}
%%            0           \\
%%            -\frac{1}{2}(u^2 + v^2) - \frac{F^{-1/2}}{\epsilon}\eta           \\
%%            \eta v                      
%%   \end{bmatrix}
%% + &\lambda_4 \begin{bmatrix}
%%            -\frac{1}{2}(u^2 + v^2) - \frac{F^{-1/2}}{\epsilon}\eta       \\
%%            0           \\
%%            \eta u                      
%%   \end{bmatrix}

Eq.~\ref{eq:RSWNullSpace} represents the subspace of perturbations to
the rotating shallow water dynamics (Eq.~\ref{eq:RSW}) that conserve
energy. The reduced coefficient space $\mathcal{P}^E$ is
18-dimensional and may be parameterized by the coordinates
$(\lambda_1 \dots \lambda_{18})$. Note that demanding energy
conservation in the dynamics has reduced the allowable dimensionality
of the LDO from $3Q=408$ (the dimensionality of $\mathcal{P}$) to 18
(the dimensionality of $\mathcal{P}^E$).

A few comments are in order. First, we note the physical significance
of the 18 vectors. The first two vectors represent scenarios in which
a perturbation is made to potential energy that is equally and
oppositely balanced by a perturbation to kinetic energy, resulting in
a net zero addition to total energy. Meanwhile, the 16 vectors
associated with $(\lambda_3 \dots \lambda_{18})$ all involve a stencil
element multiplied by $\bm{u}^{\perp}$, and hence represent rotational
terms (with zero perturbation to potential energy). We have only
defined kinetic energy as translational kinetic energy and hence
rotational terms do not affect that quantity, which explains why they
appear in Eq.~\ref{eq:RSWNullSpace}. 

Second, from a numerical perspective, many of the 16 rotational
perturbations -- while theoretically permissible -- may not make
sense. For example, consider the perturbation $[-v,u,0]u^R$ (where
$u^R$ is $u$ at the right stencil point). Numerically, there are few
reasons why a discretized PDE would perturb the right stencil value
completely independently of the other stencil values. This motivates
the possibility of introducing additional numerical constraints. For
example, we may argue that we are only interested in that subset of
the 16 rotational terms that are associated with discretized spatial
operators, such as $\nabla \bm{u}$ and $\nabla^2 \bm{u}$. This would
further reduce the dimensionality of allowable LDO coefficients. This
reduction could be achieved, for example, by switching to a less
expansive basis, such as Eq.~\ref{eq:LDOsummation2}.

\section{Reduced-Order Modeling of the LDO}
\label{sec:LDOROM}
In this section, we turn our attention to the process of constructing
a reduced-order model (ROM) for a LDO. The practical setting we
envision for this is one in which we are attempting to infer a LDO
that matches some QOI from some empirical dataset for
which direct regression is not possible, because we cannot observe the
full state. We envision that we have, however, been able to analyze
one (or several) numerical model(s) into the LDO feature space using
the constrained regression techniques above; now, we wish to search in
a neighborhood of feature space close to those numerical models for a
new LDO that produces the best match as judged by some quantity of
interest from the empirical dataset.

Because this is to be done in a Bayesian setting in which candidate
LDOs are forward integrated, we have a clear motivation to investigate
whether a ROM could reduce the computational runtime of testing a
LDO. If a ROM -- identified for a particular set of LDO coefficients
and a particular set of initial conditions -- provides reasonably
correct output over a relevant range of LDO coefficient perturbations,
then we may use that single ROM in place of the full LDO corresponding
to a different set of LDO coefficients, and perform LDO coefficient
inference on the output from the single ROM.

\new{We should point out that the ROM we use here -- POD-DEIM-Galerkin
  projection -- does not respect any conservation laws, even if the
  LDO that generated it does. This is less than optimal, given that we
  have explicitly introduced conservation statements into the
  problem. However, two remarks are in order. First, the hope is that
  the ROM will be accurate enough over the window of time during which
  the user collects data that the discrepancies between the candidate
  QOI and the true QOI will be dominated by differences in the
  structure of the physics, which is what we are ultimately attempting
  to correct. We will see some support for this hypothesis in the
  numerical examples to follow. Second, there is recent research in
  model reduction that can ensure that a ROM inherit the same
  conservation properties that were defined on the original
  model~\cite{an2008optimizing,carlberg2018conservative,carlberg2015preserving,farhat2014dimensional,farhat2015structure,karasozen2017energy,lall2003structure}. These
  structure-preserving ROMs could be used in place of the
  POD-DEIM-Galerkin approach. We leave this as an improvement to be
  explored in future research.}

We assume that we can write the relevant system dynamics for the
discretized system as follows:
\begin{equation}
\bm{\dot{X}} = \bm{f}(\bm{X}) \; ,
\label{eq:FullLDO}
\end{equation}
where $\bm{X} \in \mathbb{R}^{N_X \times N}$ is the full discrete
state and $\bm{f}(\bm{X})$ is a LDO (applied at all $N_X$ spatial
points).

One approach for building a ROM for a LDO follows the traditional
notions of POD-Galerkin projection~\cite{HolmesBook,Lumley,ChenThesis}. If we
have collected some data from simulations corresponding to a
particular LDO whose coefficients we have already identified, then we
may reduce the state dimension by computing a set of orthogonal basis
vectors from the data snapshots and projecting the state onto these
modes:
\begin{equation}
\bm{X} \approx \Phi \bm{c} \; ,
\end{equation}
and similarly, we may project the system dynamics onto these modes:
\begin{equation}
\bm{\dot{c}} \approx \Phi^T \bm{f}(\Phi \bm{c}) \; ,
\label{eq:PODGalerkin}
\end{equation}
where $\Phi \in \mathbb{R}^{N_X \times m}$, $m \ll N_X$, is a matrix whose
columns are the $m$ POD modes, and $\bm{c} \in \mathbb{R}^m$ are the
expansion coefficients.

While this process reduces the state dimension, we must still evaluate
the expensive LDO $\bm{f}(\Phi \bm{c})$ at all $N_X$ spatial points at
each timestep, so computationally we have not gained anything in going
from Eq.~\ref{eq:FullLDO} to Eq.~\ref{eq:PODGalerkin}. One strategy
for overcoming this is to form an approximating function for
$\bm{f}(\Phi \bm{c})$ based on an evaluation of it on a sparse set
$d \ll N_X$ points. The Discrete Empirical Interpolation Method
(DEIM) provides one way of achieving this~\cite{barrault2004empirical,DEIM}. Specifically,
we can approximate the LDO as follows:
\begin{equation}
\bm{f}(\Phi \bm{c}) \approx U \left( D^T U \right)^{-1} D^T \bm{f}(\Phi \bm{c}) \; ,
\end{equation}
Here, the LDO $\bm{f}(\cdot)$ is only evaluated at $d$ spatial points,
which are selected by a greedy algorithm. $U \in \mathbb{R}^{N_X
  \times d}$ is the matrix of POD modes computed for the LDO
$\bm{f(X)}$, which can be found by applying the LDO to the simulation
data. $D \in \mathbb{R}^{N_X \times d}$ is an index matrix \new{which
  corresponds to selected columns of the $N_X \times N_X$ identity
  matrix, and it is} constructed such that the product $D^T
\bm{f}(\Phi \bm{c}) = \bm{f}((\Phi \bm{c})^k)$, where $k = i_1 \dots
i_d$ are indices for the $d$ DEIM spatial points.

Altogether, the LDO ROM may be written:
\begin{equation}
\bm{\dot{c}} \approx \left[ \Phi^T U \left( D^T U \right)^{-1} \right] D^T \bm{f}(\Phi \bm{c}) = R \bm{f}_d (\bm{\Phi c}) \; ,
\end{equation}
where the ROM matrix $R = \Phi^T U \left( D^T U \right)^{-1}$ may be
precomputed, and $\bm{f}_d \equiv D^T \bm{f}$\old{ is the LDO evaluation at
the DEIM points only}.

Importantly, it should be noted that this ROM is identified for a
particular LDO, using data obtained with a particular set of initial
conditions. If either of these things \old{are} \new{is} altered, one should expect
the accuracy of the ROM to degrade. However, the ROM may still produce
output that is accurate enough for the purposes of inferring
perturbations to the LDO from which it was derived.

\section{Numerical Examples: LDO Inference}
In this section, we present a number of numerical examples
demonstrating LDO inference and structural UQ. Specifically, we use
the example of the rotating shallow water equations and the
conservation of energy law (Eq.~\ref{eq:RSWEnergy}) that they
obey. The reason we choose to use the rotating shallow water equations
is because they are relatively accessible to implement while still
preserving key ingredients of oceanic flows, such as
quasi-two-dimensionality, geostrophic and hydrostatic balance, and
topographic effects.

Our goal is \new{to} illustrate the twofold objectives of system identification
and structural UQ. Regarding the former, we wish to demonstrate how
LDO coefficients may be inferred via regression in an appropriate
feature space when the full-state is known. This is intended to
represent the practical setting in which we wish to analyze a
numerical model into the LDO basis. Regarding the latter objective, we
wish to demonstrate how the constrained LDO coefficient space
$\mathcal{P}^E$ may be searched in an area close to the numerical
model coefficients for candidate LDO coefficients that produce good
matches with some QOI from some empirical dataset.

\begin{figure}%[h]
\centering
\includegraphics[width=0.9\textwidth]{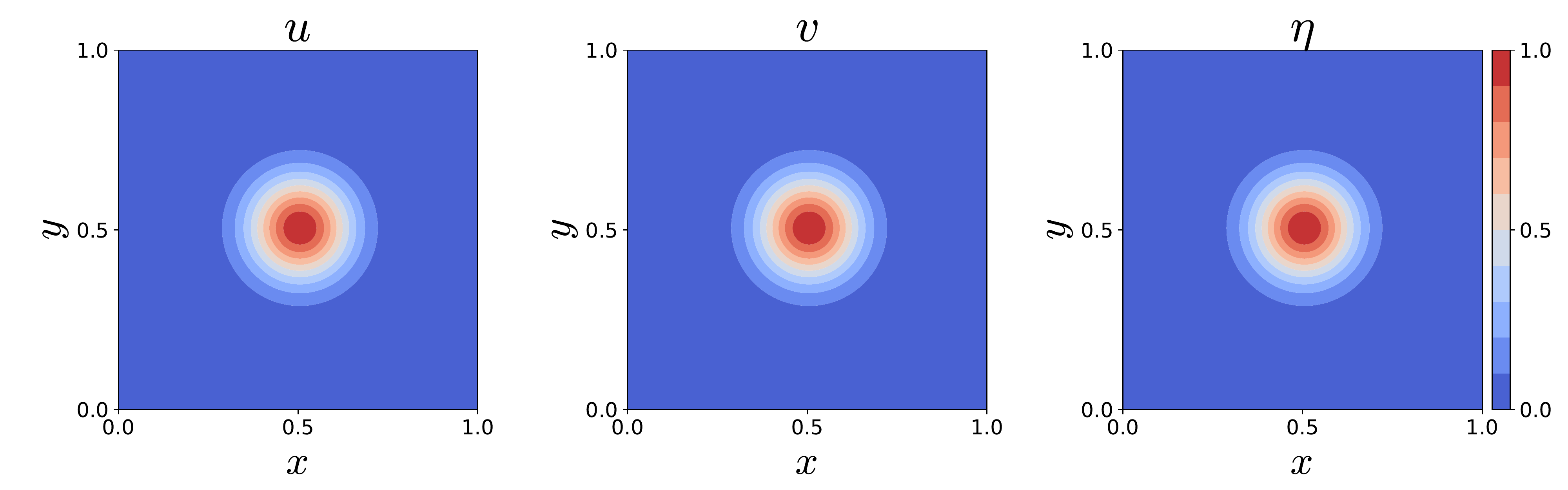}
\caption{Initial conditions used for all dynamics.}
\label{fig:InitialConditions}
\end{figure}

The LDO coefficient subspace in which we will work takes the form: 

\begin{equation}
\bm{f} = \bm{f}_{RSW} + 
  \lambda_1 \begin{bmatrix}
           -\frac{1}{2}u^2 - \frac{F^{-1/2}}{\epsilon}\eta        \\
           -\frac{1}{2}uv           \\
           u\eta                      
  \end{bmatrix}^T
+ \lambda_2 \begin{bmatrix}
           -\frac{1}{2}uv           \\
           -\frac{1}{2}v^2 - \frac{F^{-1/2}}{\epsilon}\eta           \\
           v\eta                      
  \end{bmatrix}^T \; ,
\label{eq:RSWPerturb}
\end{equation}
where $\bm{f}_{RSW}$ denotes the LDO for the rotating shallow water
dynamics (using $F = 1000$ and $\epsilon = 0.05$) and $(\lambda_1 ,
\lambda_2)$ represent perturbations in $\mathcal{P}^E$. In this
section, the rotating shallow water equations represent the numerical
model that we will learn via full-state regression, and $(\lambda_1 ,
\lambda_2) = (20,-20)$ represents the LDO corresponding to the
dynamics that underlie our empirical dataset. The initial conditions used are displayed in Fig.~\ref{fig:InitialConditions}.

We note that the 2-dimensional perturbation subspace
(Eq.~\ref{eq:RSWPerturb}) we are using (referred to hence as
$\mathcal{P}^E$) is itself a subspace of the more general
18-dimensional energy-conserving subspace derived in
\S[\ref{sec:Constraints}]: we are neglecting the 16 perturbation
vectors associated with pure rotation and are only considering those
which involve equal and opposite perturbations to kinetic and
potential energies. We have chosen to further constrain the allowable
LDO coefficient space in this way so as not to obfuscate the intended
purpose of these examples. Our goal in these examples is simply to
demonstrate LDO inference given different observables (e.g., full or
partial state observables, coarse state observables) and different
approaches to generating those observables (e.g., full LDO calculation
or ROM approximation of the LDO). In principle, consideration of the
full 18-dimensional space would change nothing about our inference
methods; it would simply require more computational time and would be
less easy to visualize.

\begin{figure}%[h!]
\begin{subfigure}[t]{0.99\textwidth}
\centering
\includegraphics[width=0.75\textwidth]{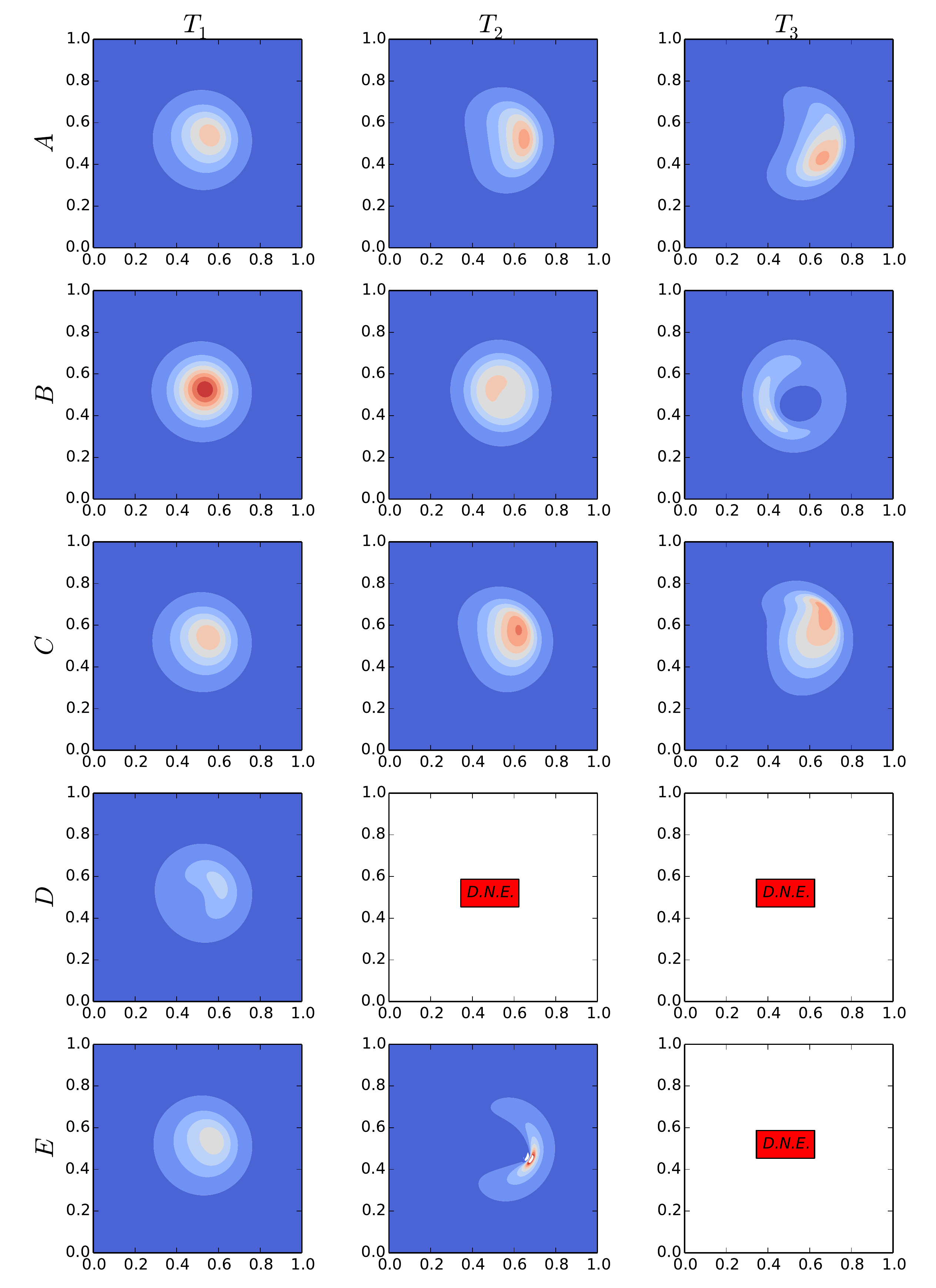}
\end{subfigure}
\begin{subfigure}[t]{0.99\textwidth}
\centering
\includegraphics[width=0.5\textwidth]{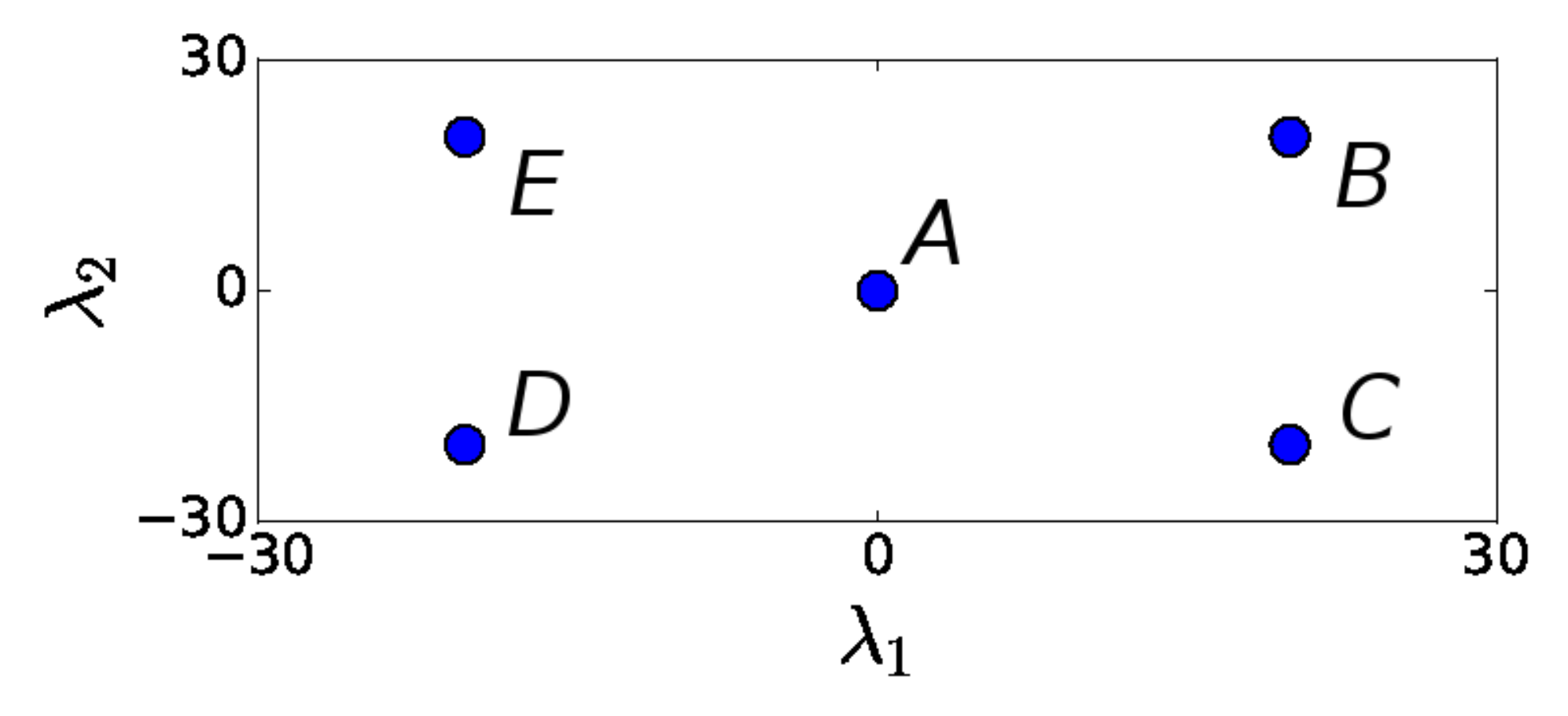}
\end{subfigure}
\caption{Illustration of effects of structural perturbations to the RSW equations on the forward-integrated dynamics. Shown are snapshots of $\eta$ for different LDOs (different values of $\bm{\lambda}$), after \new{$T_1$ = }1,000, \new{$T_2$ = }3,000, and \new{$T_3$ = }5,000 time steps with $\Delta t = 0.2\frac{1}{\Delta x^2}$. ``D.N.E.'' stands for ``Does not exist'' \new{(due to solution ``blow up'')}. The colorscale range is (-0.15,1.75).}
\label{fig:RSWvsPerturb}
\end{figure}

To illustrate the different dynamics possible in the subspace
$\mathcal{P}^E$, we display snapshots of the state variable $\eta$ for
the perturbed dynamics compared to those of the RSW equations in
Fig.~\ref{fig:RSWvsPerturb}. As
can be seen, the perturbed dynamics can look quite different relative
to the base RSW dynamics after some time has passed. Also apparent
from this is that the dynamics could potentially ``blow-up'' in
finite-time, depending on the choice of LDO coefficient
perturbations. We can get some insight into this by examination of the
energy conservation statement Eq.~\ref{eq:RSWEnergy}, which reveals
that there are terms present that are not perfect divergences and
hence can act as sources/sinks for energy. The LDO perturbations we
are allowing are chosen simply to be consistent with that differential
energy law, which is obeyed by the RSW equations. Since that
differential law does not guarantee that the total energy (i.e.,
integrated over space) will remain constant in time, we have no
guarantee that the dynamics consistent with that law should exist for
all time.

What follows are several demonstrations in which we attempt to infer
the dynamics given different observables and different models for
generating observables. Our first example is a simple demonstration of
how the LDO may be regressed from data if one has access to the full
state. After this, the remaining examples involve inferring the
perturbed dynamics under the assumption of imperfect knowledge of the
state and/or dynamics. In particular, we will investigate whether
``coarse'' (to be defined) state output is sufficient for inference
purposes, and whether a POD-DEIM ROM identified for the RSW LDO can be
used as a fast model for LDO inference.

\subsection{Full-State Data Regression}

\begin{figure}%[h!]
\centering
\includegraphics[width=0.6\textwidth]{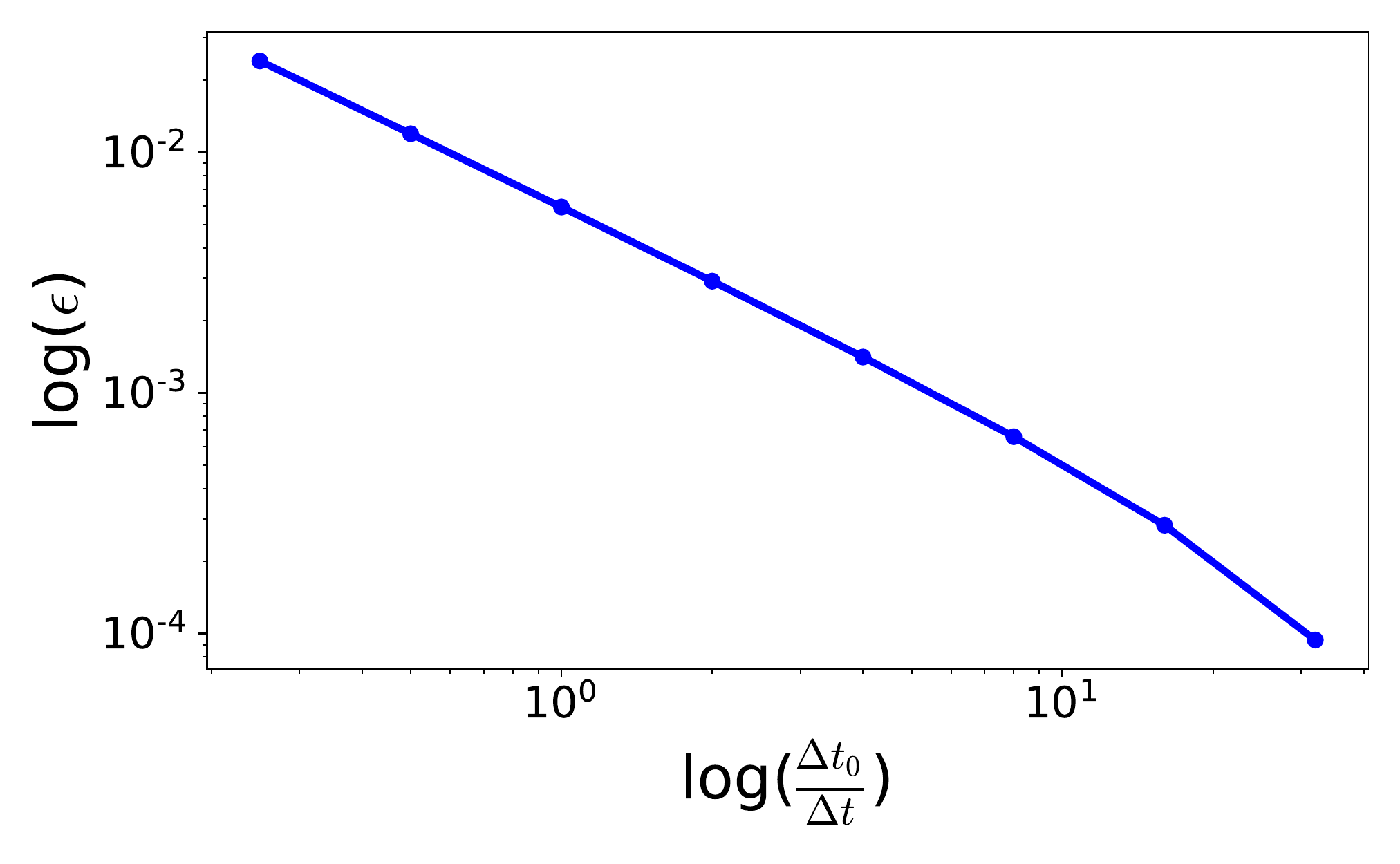}
\caption{Timestep convergence study for the forward Euler method applied to the RSW equations.}
\label{fig:convergence}
\end{figure}

Our first example is a simple confirmation that a LDO may be regressed
from data, as described in \S[\ref{sec:dataregress}]. Here, we use
data generated from the RSW equations, with no perturbation to the
dynamics. We generate data from the RSW equations using a structured
Cartesian grid of 100$\times$100 spatial points. We use explicit
timestepping with a time step $\Delta t = 0.2\frac{1}{\Delta
  x^2}$. \new{We perform a simple timestep convergence study in
  Fig.~\ref{fig:convergence}, in order to confirm that this value is
  within the asymptotic region of convergence. In this plot, $\Delta
  t_0$ corresponds to the time step value we have chosen to use. We
  forward integrate the RSW equations for a prescribed period of time
  for different values of the time step $\Delta t$, and we calculate
  the L2-norm of the solution error between a solution obtained with
  that $\Delta t$ and the same solution calculated with high temporal
  resolution (corresponding to $\frac{\Delta t_0}{\Delta t} = 64$ in
  the figure). Indeed, we observe that the L2-norm of the solution
  error decreases linearly on a log-log scale with a slope of
  approximately -1, which indicates algebraic solver convergence~\cite{wang2007high}.}

We begin by using least-squares regression with all quadratic
combinations of the stencil elements as the basis $\bm{\psi}$
\old{(see Eq. 5)}. \new{In all experiments, our stencil consists of
  the traditional 5-point stencil. For regression, we solve the normal
  equations:}
\begin{equation}
  \label{eq:normaleq}
  \new{ \left( \bm{\Psi^T \Psi} \right) \bm{P} = \left( \bm{\Psi^T Y} \right) \;\;\; , \;\;\;
  \bm{Y} = 
    \begin{bmatrix}
           \vertbar & & \vertbar \\
           \bm{f_1} & \hdots & \bm{f_N} \\
           \vertbar & & \vertbar 
         \end{bmatrix} .}
\end{equation}

\new{As described previously, $\bm{\Psi} \in \mathbb{R}^{N_S \times
    Q}$, $\bm{P} \in \mathbb{R}^{Q \times N}$, and $\bm{Y} \in
  \mathbb{R}^{N_S \times N}$. We use column-pivoted QR factorization
  to solve the problem in a well-conditioned way. Note that, in
  practice, we would begin by constraining the LDO parameter space
  with the prior physics knowledge that we introduced in
  \S[\ref{sec:Constraints}]. Indeed, we will solve the regression
  problem for this constrained case as well. However, first
  considering the unconstrained problem is a valuable pedagogical
  bridge that addresses several important issues in LDO regression
  that arise when there are many predictor variables.}

\begin{figure}%[h]
\centering
\begin{subfigure}[t]{0.48\textwidth}
\includegraphics[width=1\textwidth]{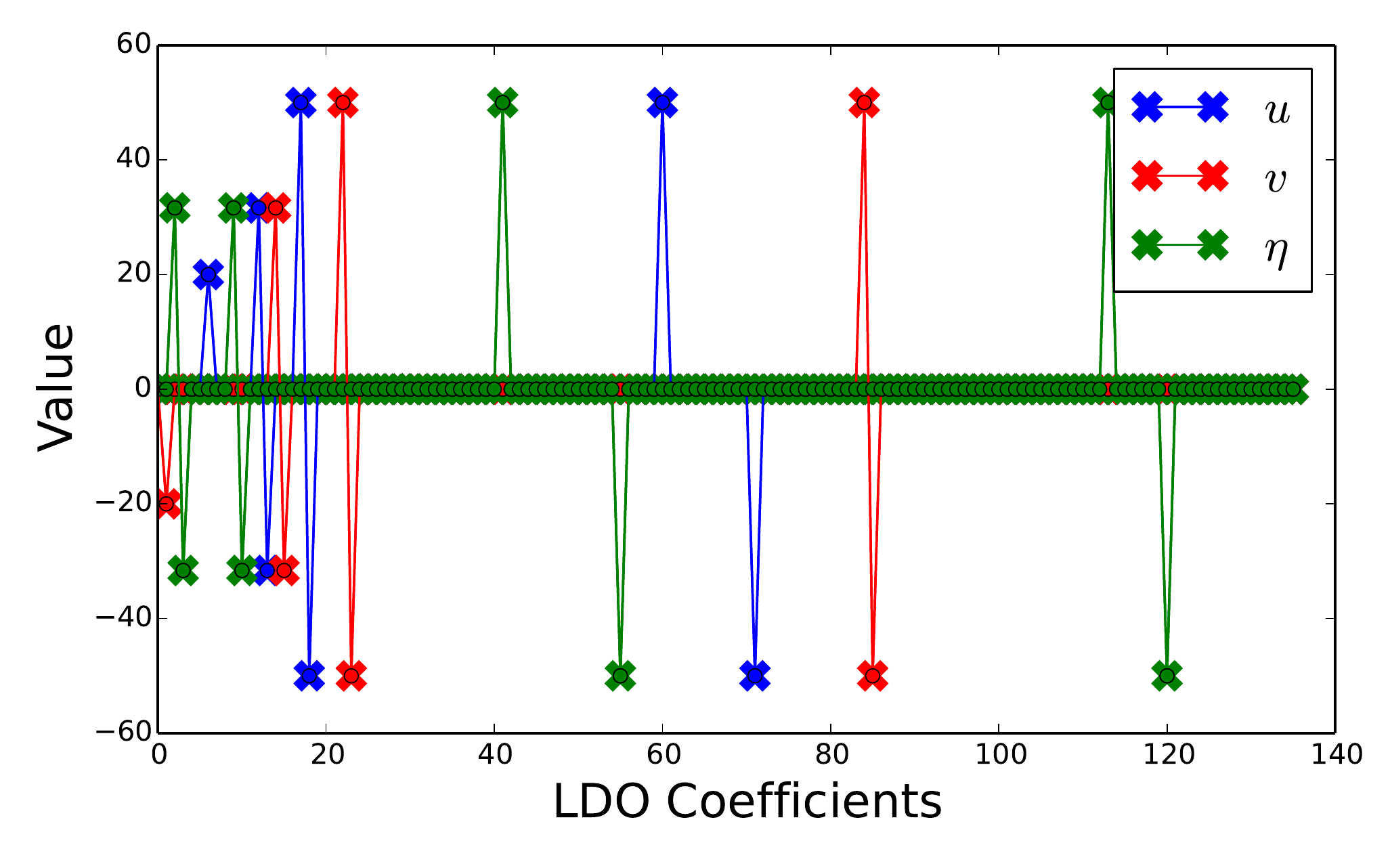}
\caption{Euler-timestepping.}
\label{fig:LDOPolyFitA}
\end{subfigure}
\begin{subfigure}[t]{0.48\textwidth}
\includegraphics[width=1\textwidth]{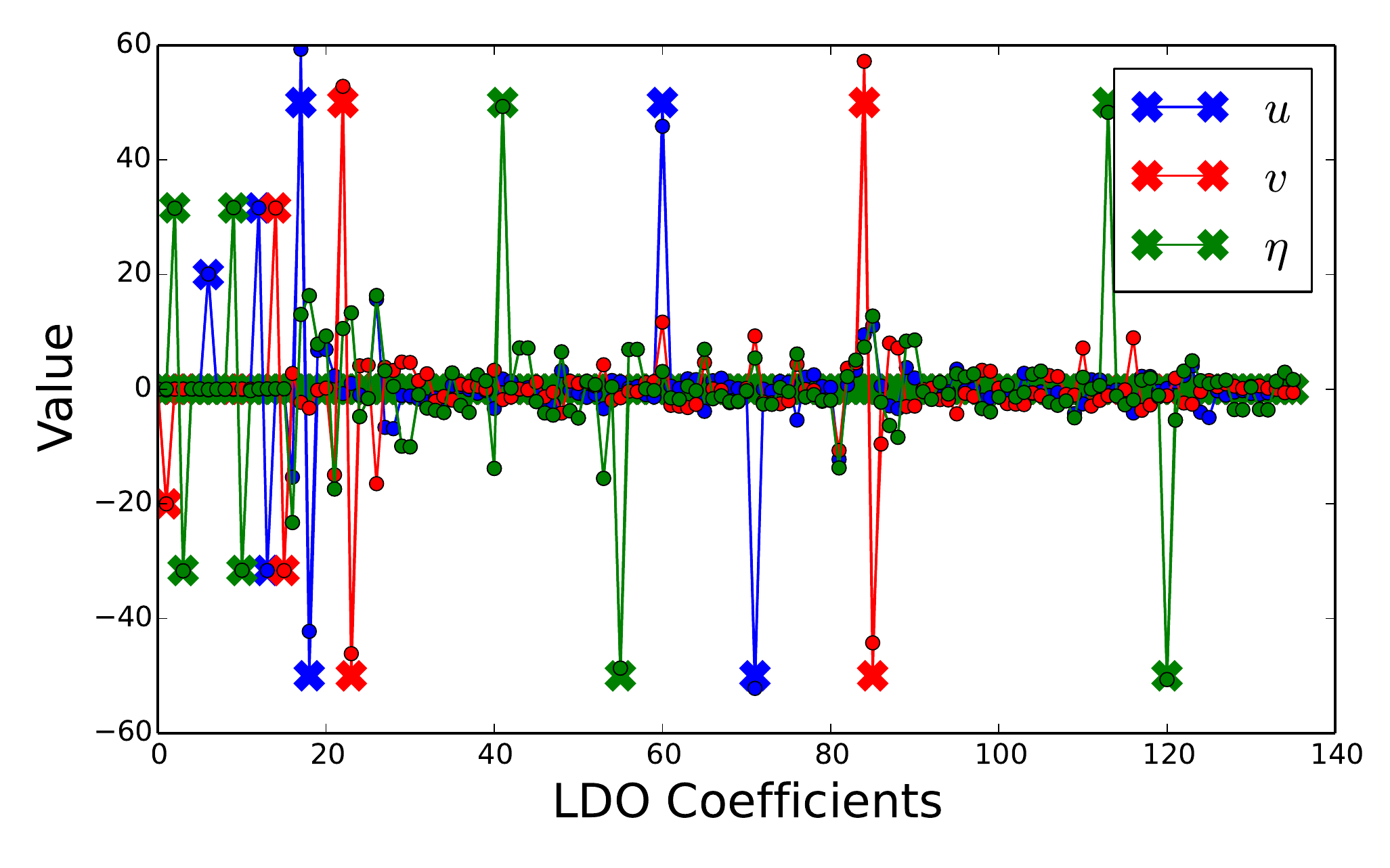}
\caption{$4^{th}$-order Runge-Kutta.}
\label{fig:LDOPolyFitB}
\end{subfigure}
\caption{Least-squares regression of the LDO coefficients for the quadratic polynomial basis ($\times$-markers denote the actual coefficients -- obtained by finite differencing the RSW equations -- and $\circ$-markers denote the regressed coefficients). \new{The basis used corresponds to Eq.~}\ref{eq:LDOsummation}.}
\label{fig:LDOPolyFit}
\end{figure}

The main difficulty encountered when \old{doing this} \new{solving the
  unconstrained regression problem} is that many of the predictors
(i.e., the basis elements) are highly-correlated with one another
(``multicollinearity''). This unfortunately makes the regression
fairly non-robust. If the data are generated using Euler
time-stepping, then the LDO should be able to exactly fit the data and
in this special case, the results are fairly good (see
Fig.~\ref{fig:LDOPolyFitA}).

However, the results are much poorer if a different numerical scheme
is used for time-stepping (e.g., Runge-Kutta). This is partially
because the LDO can no longer fit the data exactly. For example,
Runge-Kutta has the effect of both extending the stencil and enlarging
the polynomial order of the LDO, which means that no LDO consisting of
quadratic polynomial combinations of a five-point stencil will be able
to exactly fit the data. This slight error is enough to greatly
magnify the effects of the multicollinearity amongst the predictors,
resulting in a relatively poorer datafit (see
Fig.~\ref{fig:LDOPolyFitB}). It should further be noted that in the
example displayed in Fig.~\ref{fig:LDOPolyFit}, five experiments with
different initial conditions have been used in an attempt to factor
out the ``richness'' of the dataset as a contributor to
multicollinearity (indeed, if the dataset were to consist of only a
single experiment, the predictors would be more highly correlated and
the LDO coefficient fit would be even worse).

\begin{figure}%[h]
\centering
\begin{subfigure}[t]{0.48\textwidth}
\includegraphics[width=1\textwidth]{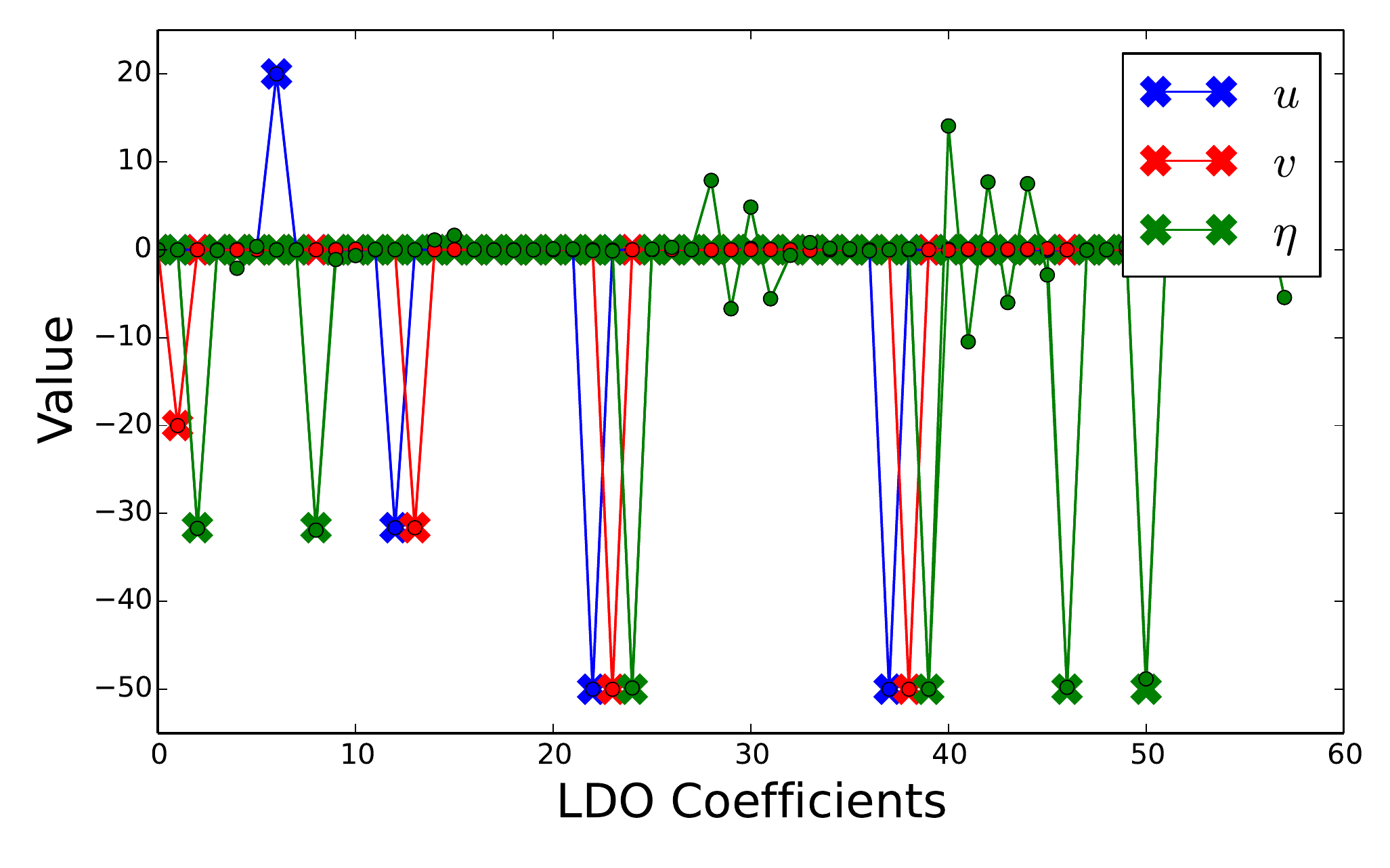}
\caption{Least squares fit.}
\end{subfigure}
\begin{subfigure}[t]{0.48\textwidth}
\includegraphics[width=1\textwidth]{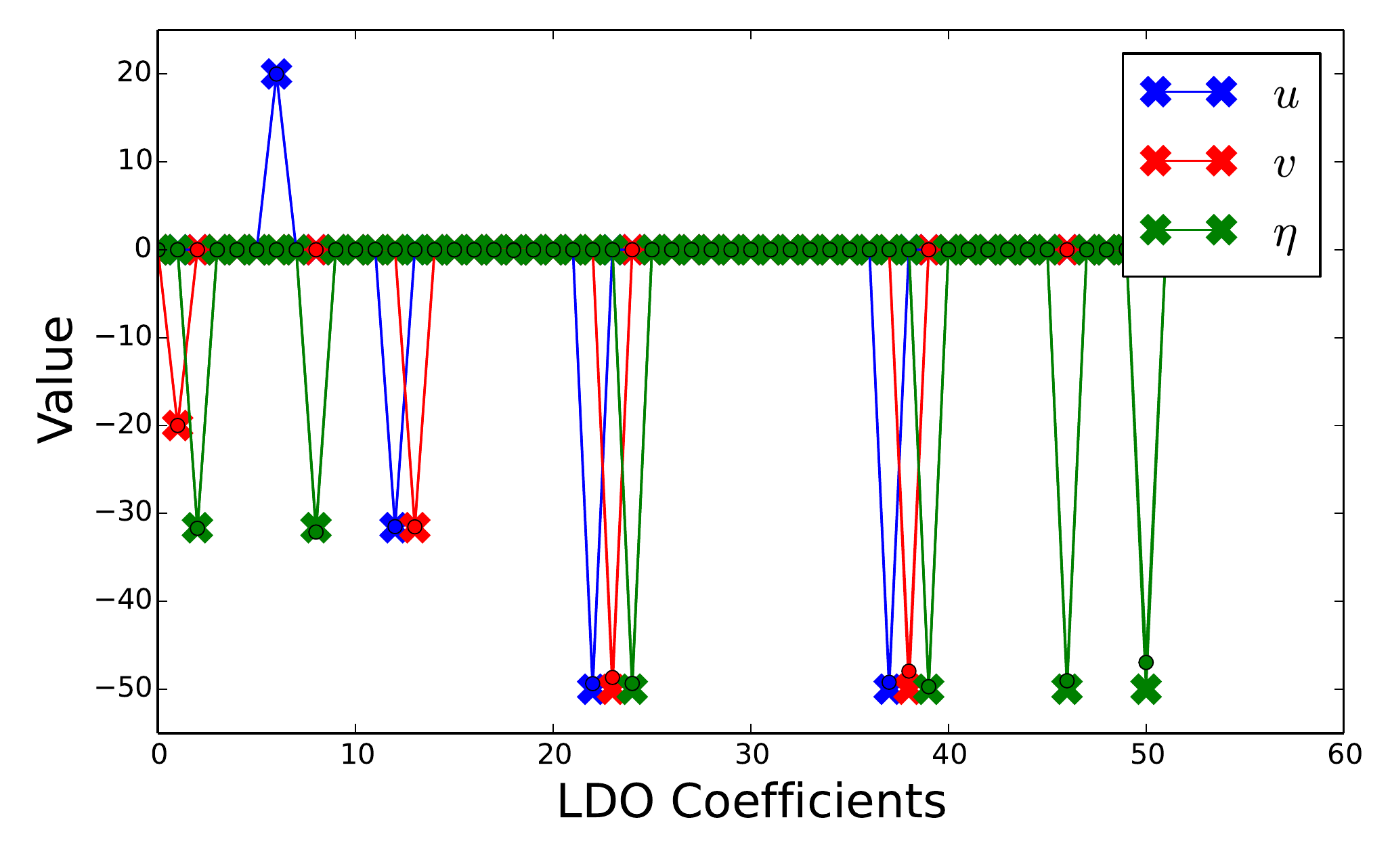}
\caption{LASSO fit.}
\end{subfigure}
\caption{LDO coefficient regression using the differential operator basis ($\times$-markers denote the actual coefficients -- obtained by finite differencing the RSW equations -- and $\circ$-markers denote the regressed coefficients). Data generated using $4^{th}$-order Runge-Kutta. \new{The basis used corresponds to Eq.~}\ref{eq:LDOsummation2}.}
\label{fig:LDOOperFit}
\end{figure}

If we hypothesize that the main hurdle to achieving a robust numerical
regression in our case is multicollinearity, then there are several
``standard'' strategies we may exploit to alleviate this issue. First,
we can switch to using a basis consisting of only some differential
operators (see Eq.~\ref{eq:LDOsummation2}); this illustrates another
important reason for using such a basis in addition to LDO basis
dimension reduction. Second, we can use the LASSO method instead of
least-squares fitting. This method promotes a sparse solution to the
LDO inference problem\old{.}\new{, which is justified on the
  hypothesis that the solution will be sparse. Note that, in practice,
  one would generally need to use a technique such as cross-validation
  to arrive at an acceptable weight parameter for LASSO.}
Fig.~\ref{fig:LDOOperFit} shows an example of how using the
differential operator basis -- along with LASSO regression -- can
greatly improve the accuracy of the inferred LDO coefficients.

\new{Ultimately, however, we are interested in solving the constrained
  regression problem, since that solution will inform the prior we use
  in the subsequent structural UQ step. The setup for this problem is
  a slight variation on what we have already done. As before, we solve
  the normal equations as in Eq.~\ref{eq:normaleq}; however, the
  matrices we use are now as follows:}
\begin{equation}
  \new{
  \bm{\Psi} = 
  \begin{bmatrix}
    \bm{\Psi}_u \\
    \bm{\Psi}_v \\
    \bm{\Psi}_{\eta} 
  \end{bmatrix} \in \mathbb{R}^{3N_S \times 2} \;\;\; , \;\;\;
  \bm{P} = \bm{\lambda} \in \mathbb{R}^{2 \times 1} \;\;\; , \;\;\;
  \bm{Y} = 
  \begin{bmatrix}
    \bm{\Delta f}_u \\
    \bm{\Delta f}_v \\
    \bm{\Delta f}_{\eta}
  \end{bmatrix} \in \mathbb{R}^{3N_S \times 1} .}
\end{equation}

\new{Here, the submatrices $\bm{\Psi}_u, \bm{\Psi}_v,\bm{\Psi}_{\eta}
  \in \mathbb{R}^{N_S \times 2}$ represent the $(u,v,\eta)$ components
  of the constrained basis vectors displayed in
  Eq.~\ref{eq:RSWPerturb}, and $\bm{\Delta f = f-f_{RSW}}$ (the
  subscripts denote the respective state components). For a given
  snapshot of data, $\bm{f_{RSW}}$ is calculated simply by applying
  the LDO corresponding to the RSW equations to that
  snapshot. Formulating the regression problem in this way, we
  calculate a numerical result of $\bm{\lambda} = (-0.015, -0.030)$,
  which is very close to the true values of $\bm{\lambda} = (0, 0)$.}

\subsection{Single State Observable}

We now turn our attention to the more difficult task of inferring the
perturbed LDO given only imperfect knowledge of the state. This
problem is meant to illustrate structural UQ -- it is a scenario in
which we only have access to incomplete information from an empirical
dataset, and we wish to find a LDO that matches that QOI. The
full-state inference of the RSW LDO in the previous section represents
the numerical code model that we have learned, and now our goal is to
search in an area of the energy-constrained subspace close to that
numerical model for a LDO that produces a good match with the
empirical QOI. \new{Thus, the link between the previous step
  (regression) and this one (structural UQ) comes via the prior
  distribution we impose over the constrained LDO parameter
  space. Specifically, we choose a prior distribution that is centered
  at the mean values of $\bm{\lambda}$ that we inferred last
  section. Intuitively, our prior should weight points closer to that
  mean value more highly, and so we use a Gaussian distribution with a
  standard deviation of $30$.}

The incomplete QOI that we consider in this case will be the state
variable $\eta$, and nothing else. We assume that we are able to
observe snapshots of this state variable at any point in time from the
black-box dynamics of interest, and we wish to infer the LDO from this.

We elect to use the Metropolis MCMC algorithm to perform inference of
the parameters $(\lambda_1 , \lambda_2)$.  With access to the
fully-resolved variable $\eta$, we define our observable for the
Metropolis algorithm as the RMS error of $\eta$:
\begin{equation}
\epsilon_{l_2} = \left[ \int_{0}^T \left( \int_{V} \left| \left(\eta^{\text{true}} - \eta^{\text{model}} \right)(\bm{x},t) \right|^2 d\bm{x} \right) dt \right]^{1/2} \; ,
\label{eq:MetropolisError}
\end{equation}
where $\eta^{\text{model}}$ is $\eta$ calculated from some LDO under
consideration. We approximate the integral in
Eq.~\ref{eq:MetropolisError} using 5,000 snapshots in time, collected
using a time step $\Delta t = 0.2\frac{1}{\Delta x^2}$, where $\Delta
x$ denotes the grid point spacing (we use a 100$\times$100 grid of
uniformly distributed points).

We define the Bayesian likelihood function as a normal distribution
function of $\epsilon_{l_2}$:
\begin{equation}
\pi (\bm{\lambda} | \epsilon_{l_2}) = \text{exp} \left(-\frac{1}{2} \frac{(\epsilon_{l_2})^2}{\sigma^2} \right) \; .
\label{eq:Likelihood}
\end{equation}
Of course, the parameter $\sigma$ is a free-parameter that must be
chosen; we set $\sigma = 300$. \old{We use a uniform prior in the Metropolis
algorithm.}

\begin{figure}%[h]
\centering
\includegraphics[width=0.7\textwidth]{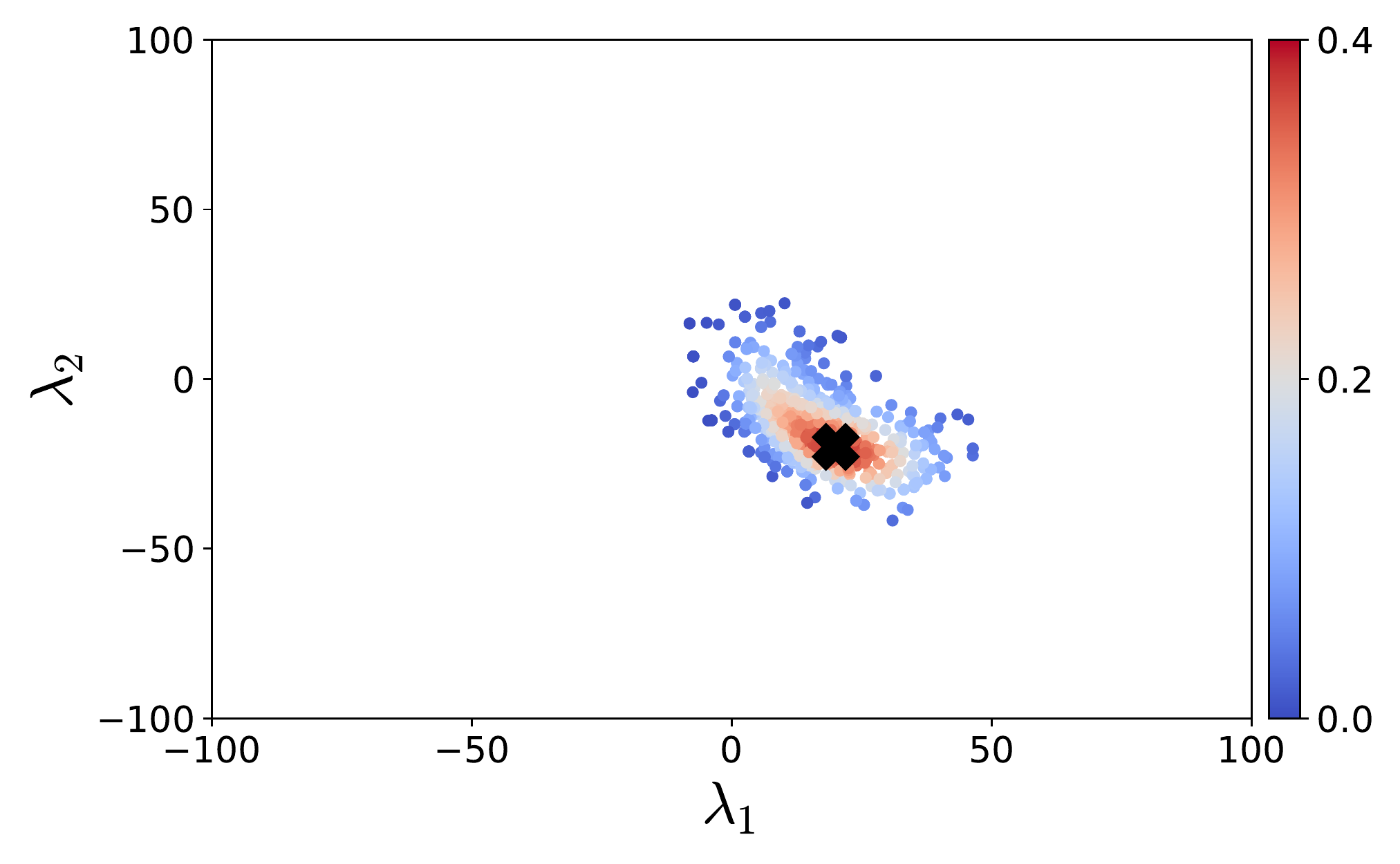}
\caption{Posterior distribution using the fully-resolved state observable for $\bm{\lambda}$ generated with 1,000 MCMC samples ($\bm{\lambda}^{\text{true}}$ is shown as the black x).}
\label{fig:2DInferenceExample}
\end{figure}

The posterior, estimated using 1,000 samples, is displayed in
Fig.~\ref{fig:2DInferenceExample}. We can clearly see that the MCMC
samples converge strongly to a cluster which is exactly centered about
the true values of $\bm{\lambda}$.

\subsection{Coarse-Data Observable}

The previous examples illustrate how constrained inference may proceed
in situations where access to one or more fully-resolved states is
assumed. This is useful for demonstration purposes of the MCMC
algorithm; however, the more realistic scenario for constrained
inference occurs when we only have access to some observable function
of the state which is of lower informational value. We show here an
example of LDO inference in the energy-constrained subspace using
state observations that are coarse in space and time. This emulates a
situation in which only sparse observations of low spatial quality are
available for structural UQ purposes.

\begin{figure}%[h]
\centering
\begin{subfigure}[t]{0.45\textwidth}
\includegraphics[width=1\textwidth]{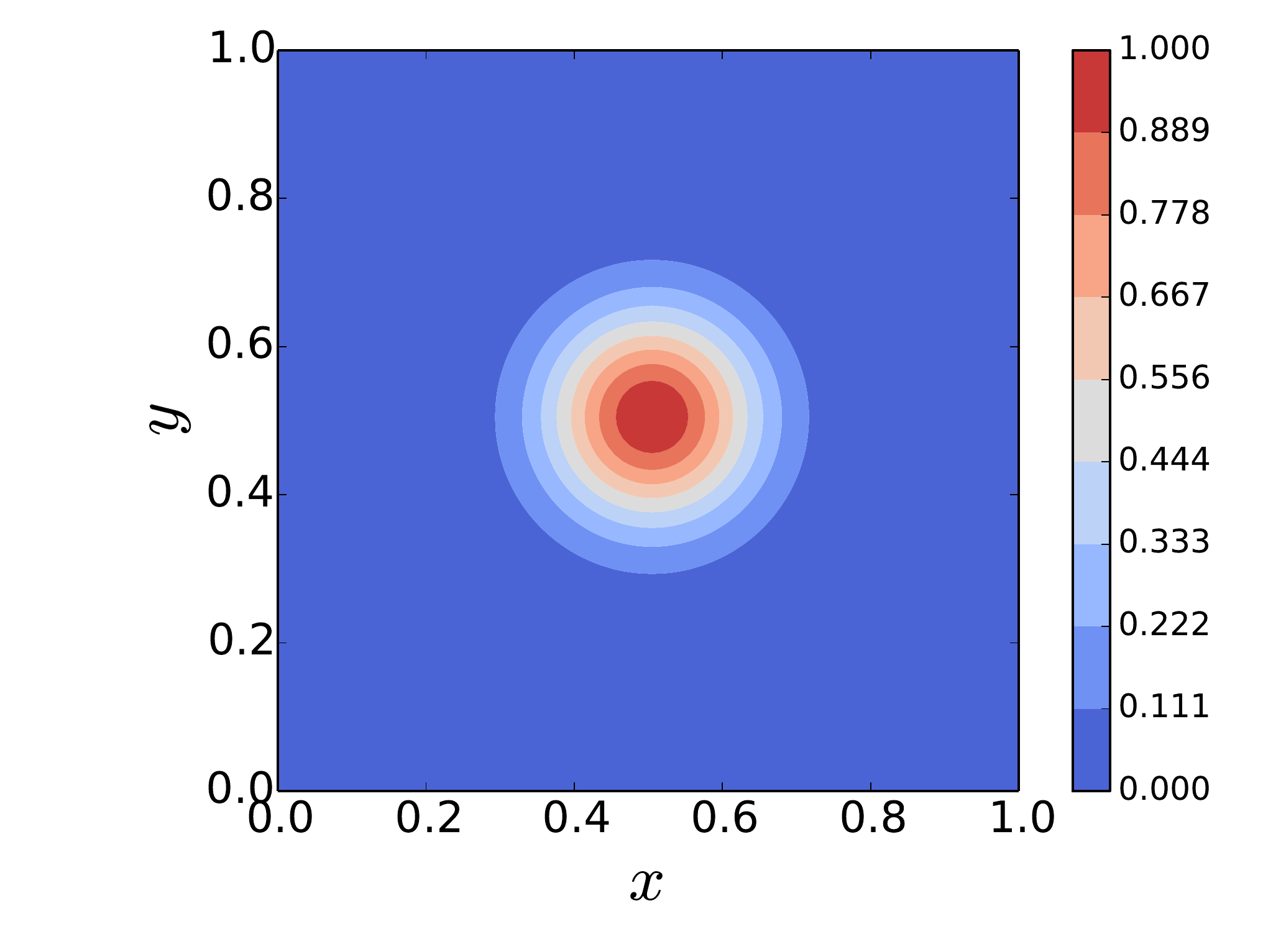}
\caption{Fine data field.}
\end{subfigure}
\begin{subfigure}[t]{0.45\textwidth}
\includegraphics[width=1\textwidth]{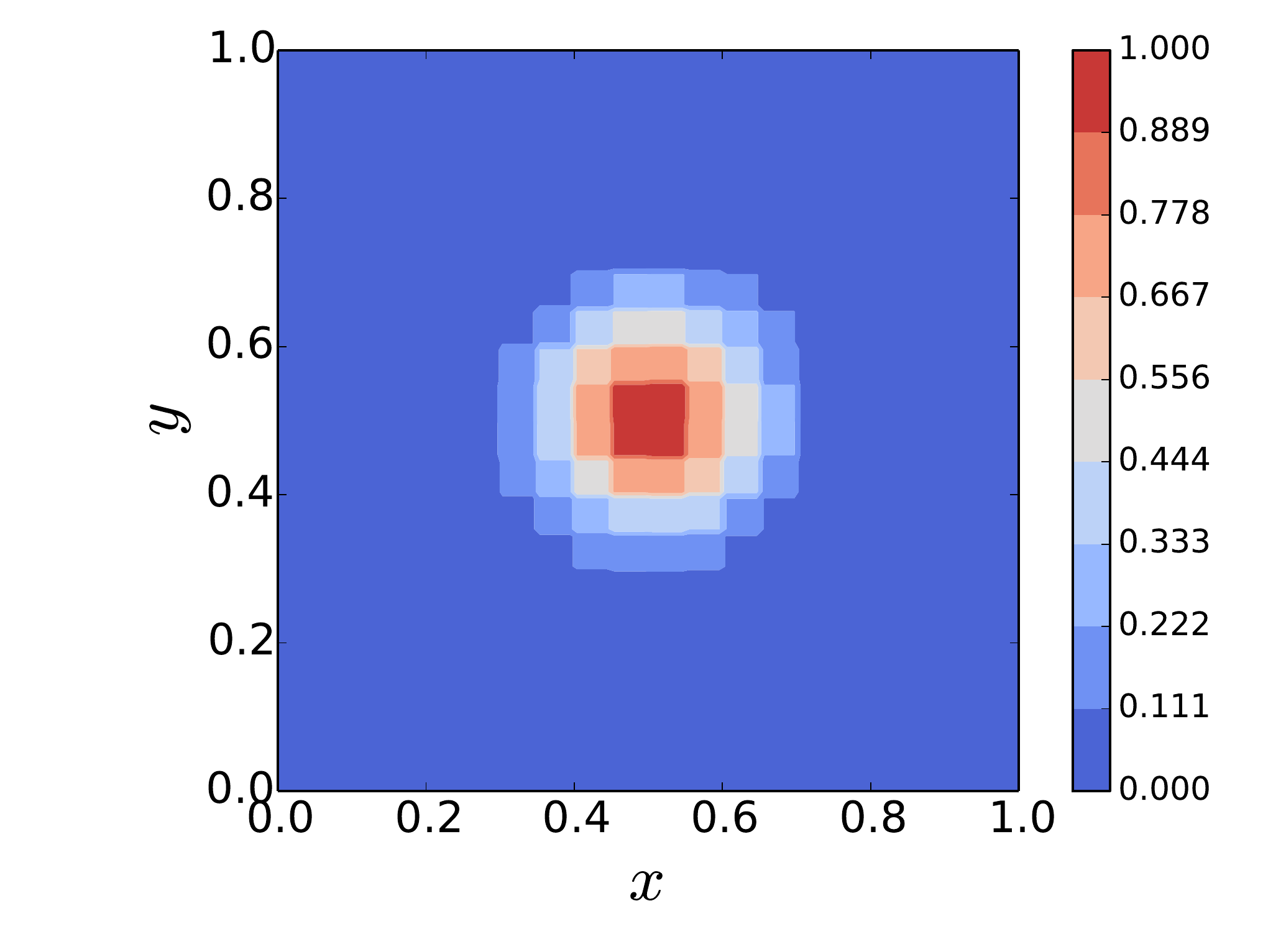}
\caption{Coarse data field.}
\end{subfigure}
\caption{Illustration of the coarsening operation.}
\label{fig:LocalAvg}
\end{figure}

In this example, we define ``coarse'' data by interpreting it as a
kernel approximation to the actual state; for example,
\begin{equation}
\eta(\bm{x'},t') \approx \mathcal{C}(\eta(\bm{x'},t')) = \int_T \int_{X} K(\bm{x'},\bm{x},t',t) \eta(\bm{x},t) d\bm{x} dt
\label{eq:KernelApprox}
\end{equation}
is an approximation to the state $\eta$. Different choices of the
kernel $K(\bm{x'},\bm{x},t',t)$ are possible. For simplicity, we use
the following discontinuous kernel:
\begin{equation}
K(\bm{x'},\bm{x},t',t) = 
\begin{cases}
    \frac{1}{|V|}, & \text{if } (\bm{x'},t') \in V \text{and } (\bm{x},t) \in V  \\
    0,                  & \text{otherwise} \; ,
\end{cases}
\end{equation}
where $V$ is a hypercube in space-time. This choice of kernel
corresponds to a local space-time average in some $V$. Thus, we can
then coarsen the data by dividing space-time into several
non-overlapping boxes, and then computing the local average within
each of these boxes via Eq.~\ref{eq:KernelApprox}. If we take this
locally-averaged data to be our observable, we may define an error
metric just as in Eq.~\ref{eq:MetropolisError}, except with
$\eta^{\text{true}}$ and $\eta^{\text{model}}$ replaced by their
coarsened versions. One way to express the level of coarsening is with
the ratio of the fine-to-coarse grid size. In this example, we use a
spatial coarsening ratio of 25 (i.e., 25 fine grid points are
contained in each coarse grid box), and a temporal coarsening ratio of
25 (i.e., 25 snapshots in time are averaged together in the coarsening
operation). Example snapshots using this process are shown in
Fig.~\ref{fig:LocalAvg}.
%% \begin{equation}
%% \begin{aligned}
%% \mathcal{C}(\eta(\bm{x'},t')) &= \int_T \int_{X} K(\bm{x'},\bm{x},t',t) \eta(\bm{x},t) d\bm{x} dt \\
%% &= \frac{1}{|V|} \int_T \int_X \mathbbm{1}_{V}(\bm{x'},\bm{x},t',t) \eta(\bm{x},t) d\bm{x}dt \\
%% &= 
%% \begin{cases}
%%     \frac{1}{|V|} \int_{V} \eta(\bm{x},t) dV , & \text{if } (\bm{x'},t') \in V \\
%%     0,                  & \text{otherwise}
%% \end{cases}
%% \end{aligned}
%% \end{equation}

\begin{figure}%[h]
\centering
\includegraphics[width=0.7\textwidth]{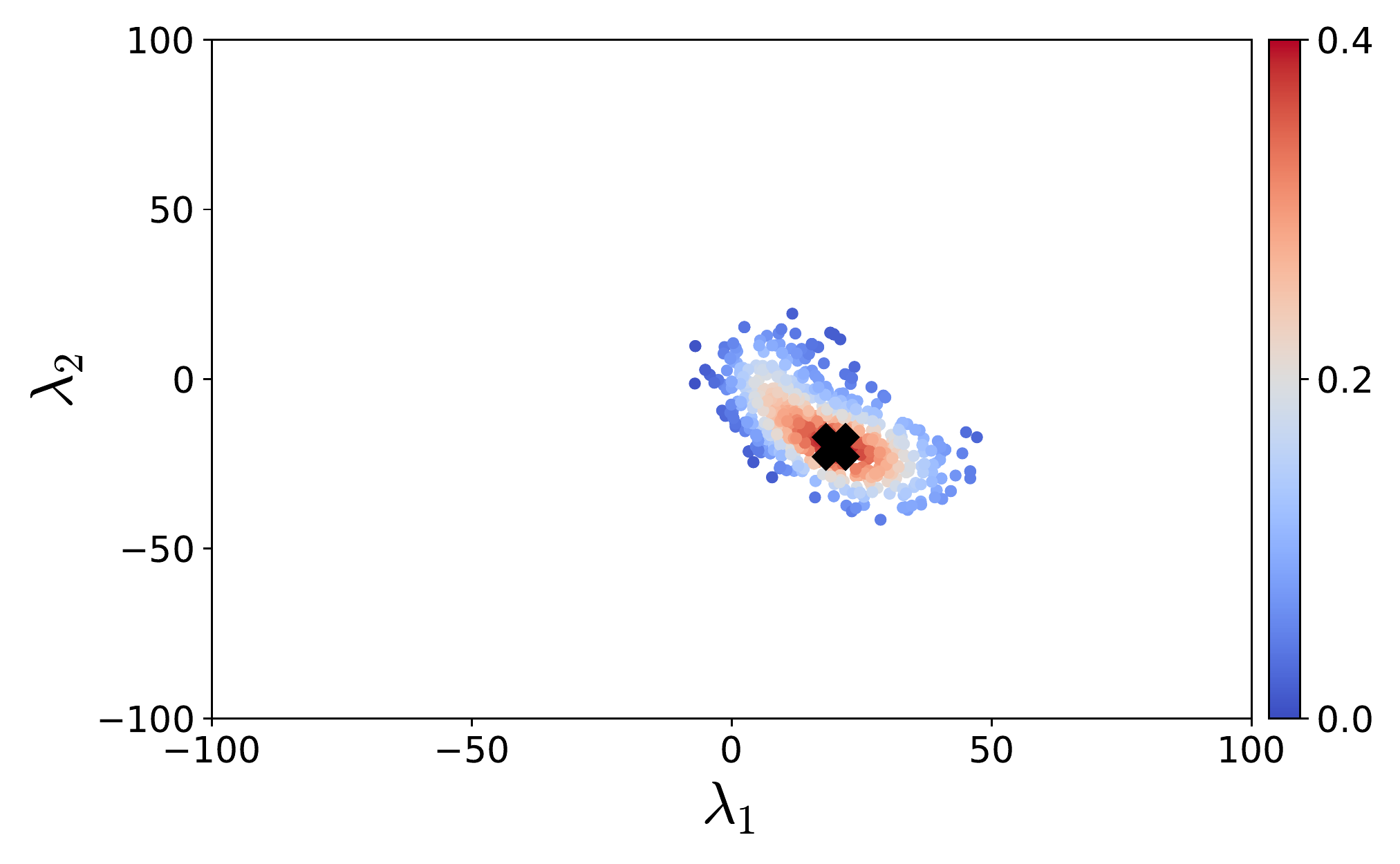}
\caption{Posterior distribution using the coarsened $\eta$ observable for $\bm{\lambda}$ generated with 1,000 MCMC samples ($\bm{\lambda}^{\text{true}}$ is shown as the black x).}
\label{fig:CoarsePosterior}
\end{figure}

As in the previous example, we perform our inference in the
energy-constrained subspace $\mathcal{P}^E$, and we are attempting to
infer the parameters $(\lambda_1 , \lambda_2)$. The MCMC posterior
(calculated with 1,000 samples, and $\sigma = 300$ in
Eq.~\ref{eq:Likelihood}) using the coarse-data observable is shown
in Fig.~\ref{fig:CoarsePosterior}. The MCMC samples converge as before
to a cluster in parameter space, which is roughly centered about the
true values.

\subsection{Coarse-Data Observable, Computed with a ROM}
In this example, we investigate whether LDO coefficient inference can
be successfully done using a ROM identified for the base LDO in place
of the actual LDO corresponding to a given perturbation. The
methodology for calculating a LDO ROM is sketched in
\S[\ref{sec:LDOROM}]. Again, the motivation for this is computational
expense: doing Bayesian inference in LDO space could be expensive,
even with energy constraints, simply because it involves calculating
discrete-time dynamics by forward-integrating candidate LDOs.

As before, our base LDO (and the ROM identified for it) corresponds to
the RSW equations (which corresponds to $\bm{\lambda} = (0,0)$). We
use the same energy-preserving perturbation to the base LDO as in the
previous example, and perform our inference in the corresponding
2-dimensional space $\mathcal{P}^E$. We take our observable to be the
state variable $\eta$, coarsened with the operator
$\mathcal{C}(\cdot)$ as defined in Eq.~\ref{eq:KernelApprox}. We use
the error metric in Eq.~\ref{eq:MetropolisError} (albeit with
the coarsened versions of $\eta$).

\begin{figure}%[h]
\centering
\begin{subfigure}[t]{0.95\textwidth}
\includegraphics[width=1\textwidth]{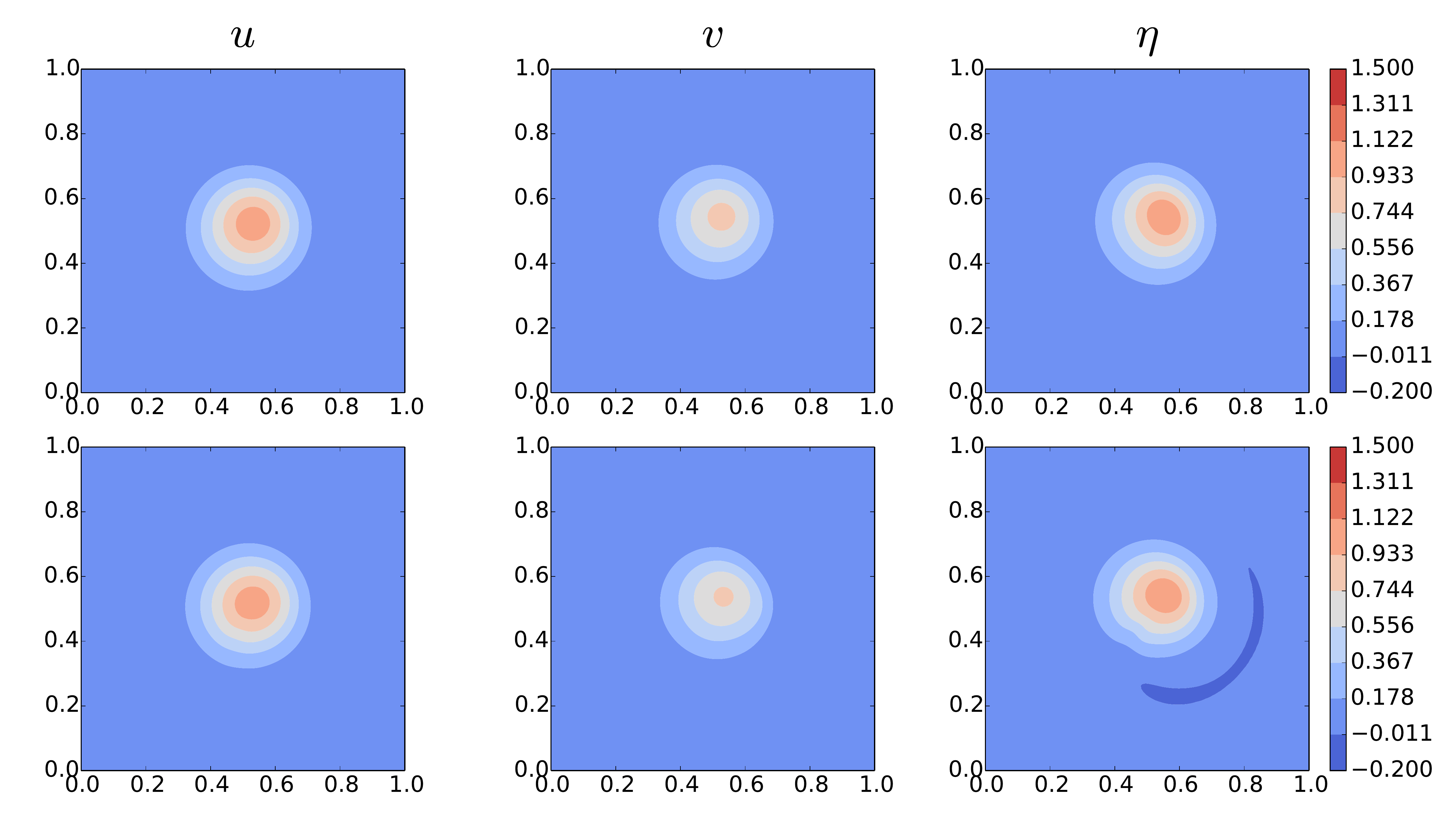}
\caption{1,000 timesteps.}
\label{fig:ActualVsROM1000}
\end{subfigure} \\
\begin{subfigure}[t]{0.95\textwidth}
\includegraphics[width=1\textwidth]{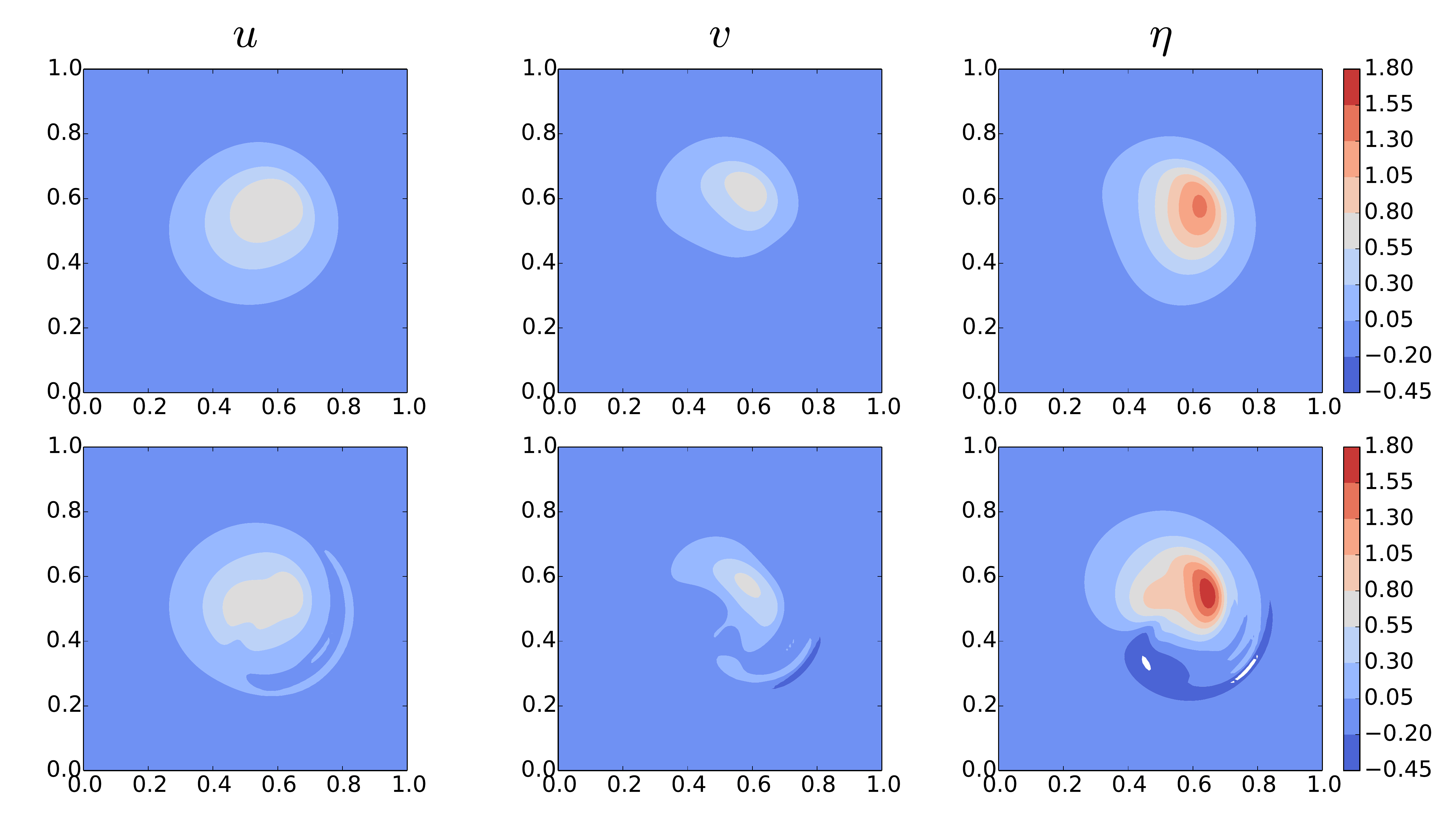}
\caption{3,000 timesteps.}
\label{fig:ActualVsROM3000}
\end{subfigure}
\caption{Snapshots of the true LDO states (rows 1 and 3) compared to ROM-calculated states (rows 2 and 4) using $\bm{\lambda} = (20,-20)$ and $\Delta t = 0.2\frac{1}{\Delta x^2}$.}
\label{fig:ActualVsROM}
\end{figure}

\begin{figure}%[h]
\centering
\includegraphics[width=0.7\textwidth]{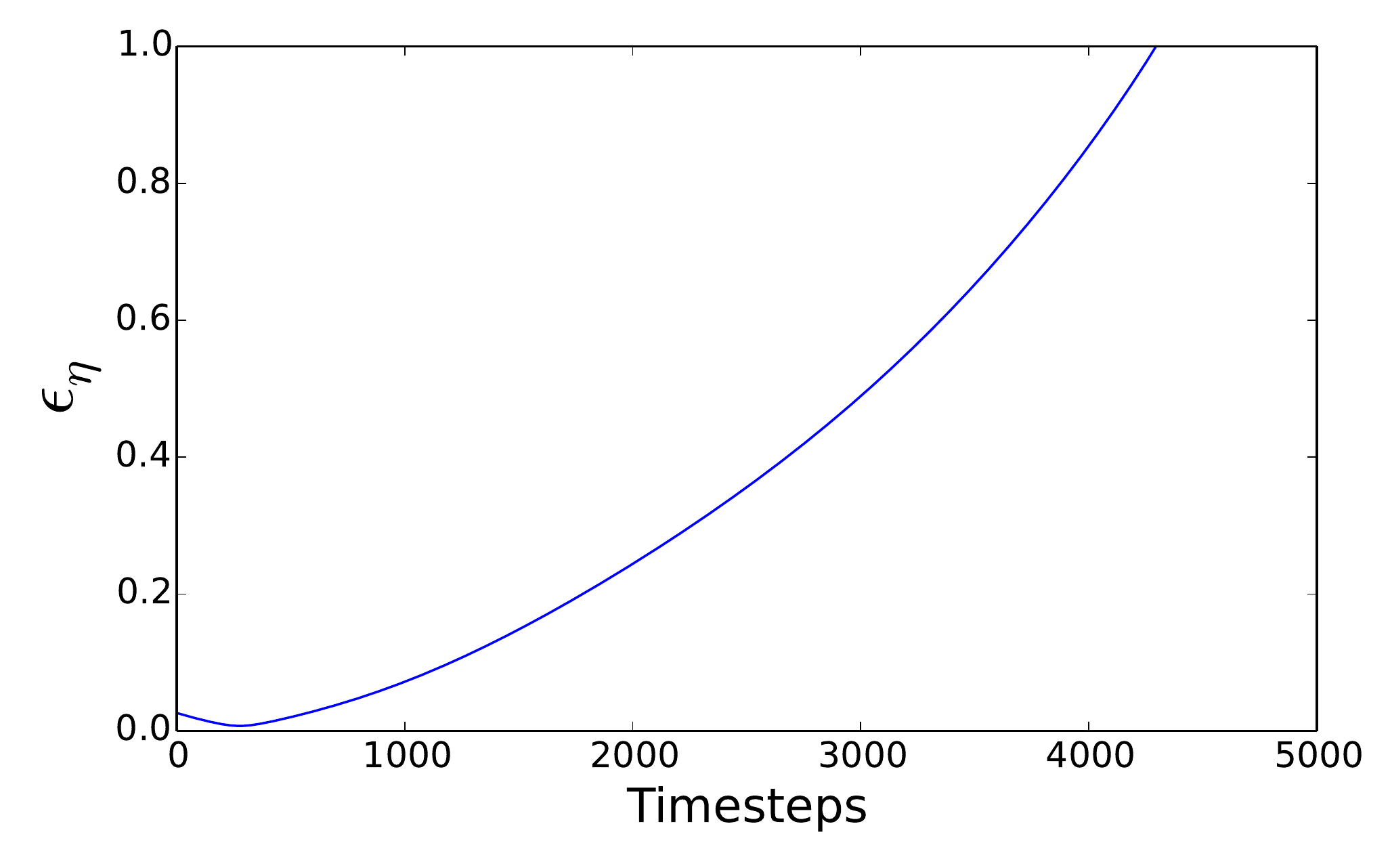}
\caption{Discrepancy between true $\eta$ ($\bm{\lambda} = (20,-20)$) and ROM-calculated $\eta$ as a function of time ($\Delta t = 0.2\frac{1}{\Delta x^2}$). Here, we use the error metric $\epsilon_{\eta} (t) = \int_{\bm{x}}|\eta_{true} - \eta_{ROM}|/|\eta_{true}| (\bm{x},t) d\bm{x}$.}
\label{fig:ActualVsROMError}
\end{figure}

To begin, we might wonder about \old{the effect of} collecting increasing
numbers of snapshots \new{from candidate LDO simulations in the MCMC algorithm, and how this affects} \old{on} the accuracy and precision of the calculated
posterior. Usually, we would expect this relationship to be monotonic:
more data is always better. However, because we are using a ROM whose
accuracy is limited by the \new{original base} data used to derive it \new{(i.e., snapshots from the RSW LDO)}, this is not
necessarily the case. That is, there is a temporal trade-off between
aggregating more information and the degrading accuracy of the
ROM. Since all simulations begin with the same initial conditions as
those used to derive the RSW ROM (shown in
Fig.~\ref{fig:InitialConditions}), the ROM is most accurate for a
given perturbed LDO at early times, when the spatial characteristics
of the state variables are close to those of the RSW equations and
hence largely spanned by the POD bases used to describe them in the
ROM (see Fig.~\ref{fig:ActualVsROM1000}). However, as time progresses,
the spatial patterns of the solutions quickly diverge from those
exhibited by the RSW equations as the differences in structural
governing dynamics become increasingly apparent. While this
information is useful in discriminating amongst different LDOs, it is
only useful if the ROM can reproduce it. As would be expected, the
accuracy of the ROM does degrade as the spatial patterns exhibited by
a given LDO diverge from those used by the ROM to describe them (see
Fig.~\ref{fig:ActualVsROM1000}, \ref{fig:ActualVsROM3000} and \ref{fig:ActualVsROMError}).

\begin{figure}%[h]
\centering
\begin{subfigure}[t]{0.47\textwidth}
\includegraphics[width=1\textwidth]{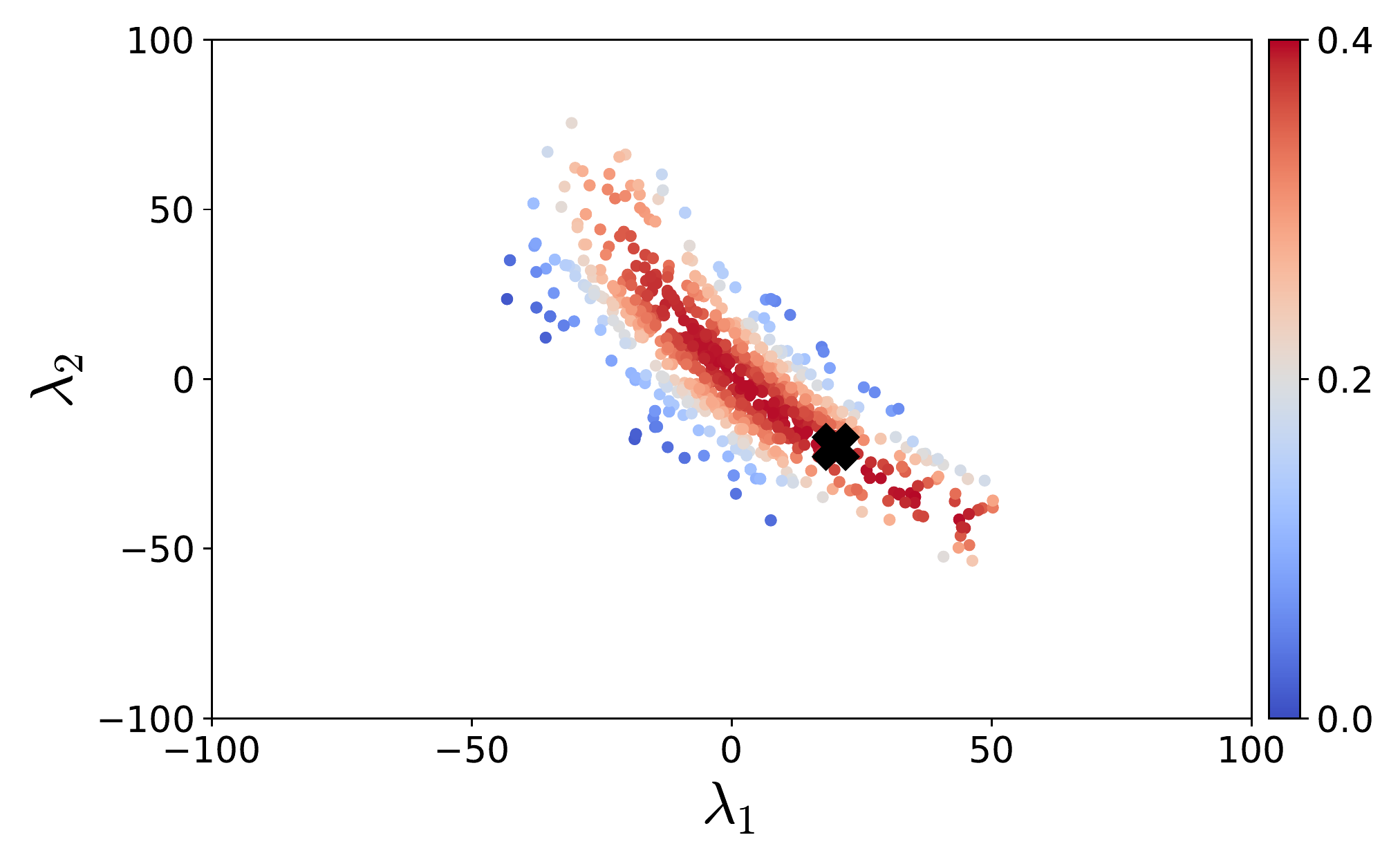}
\caption{1,000 time steps, $\sigma = 100$}
\label{fig:ROMinference1a}
\end{subfigure}
\begin{subfigure}[t]{0.47\textwidth}
\includegraphics[width=1\textwidth]{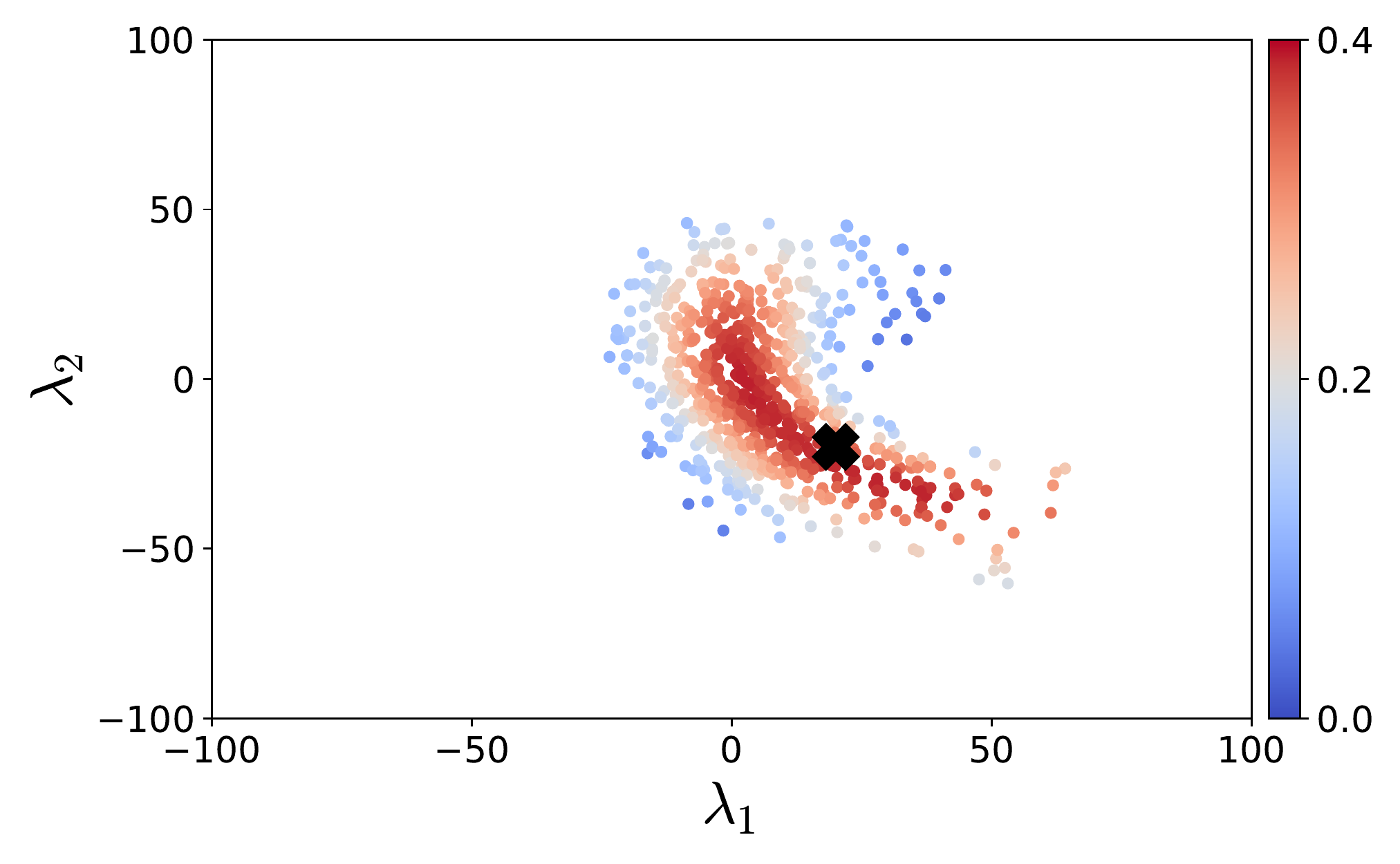}
\caption{3,000 time steps, $\sigma = 600$}
\label{fig:ROMinference1b}
\end{subfigure} \\
\begin{subfigure}[t]{0.47\textwidth}
\includegraphics[width=1\textwidth]{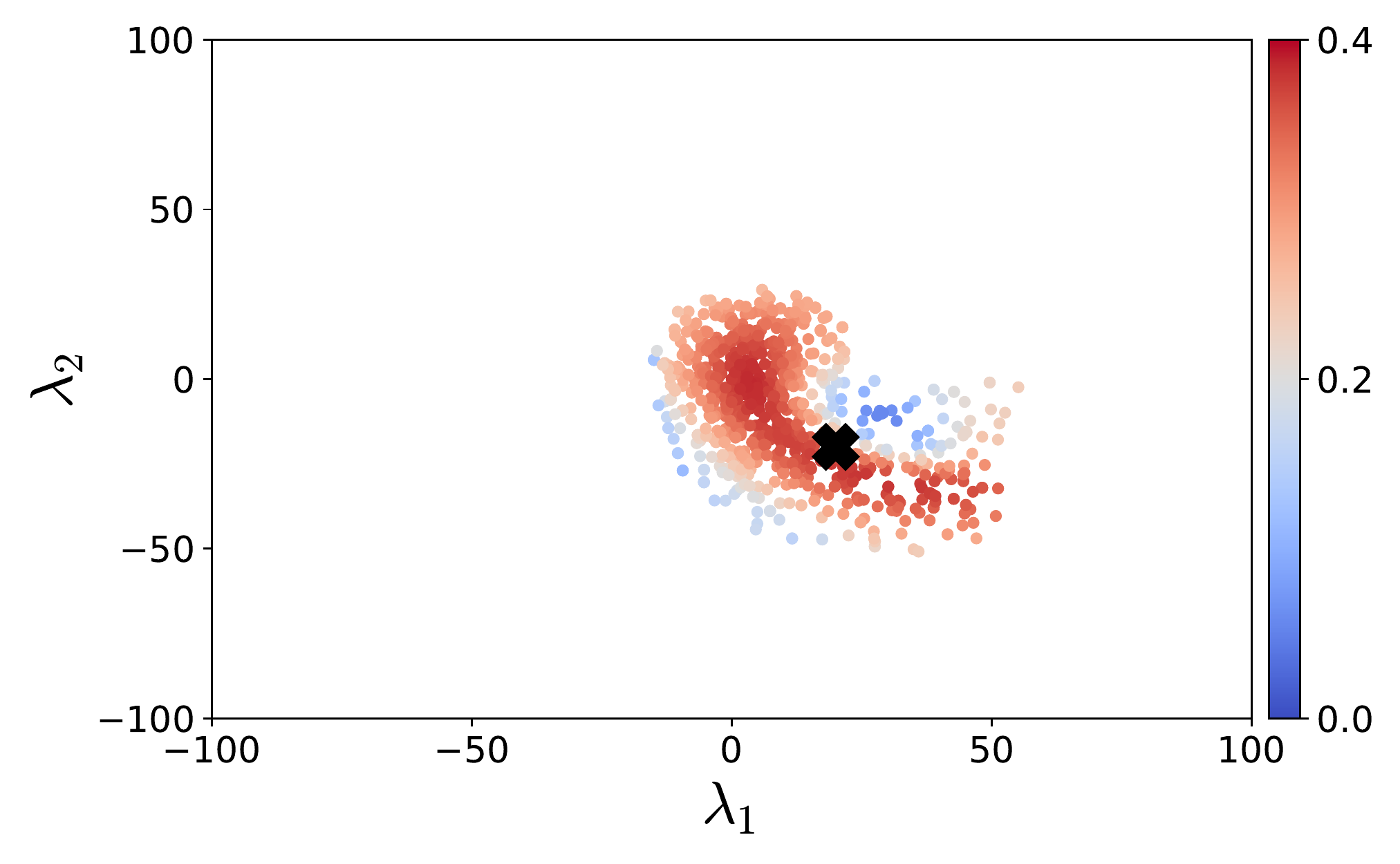}
\caption{5,000 time steps, $\sigma = 1500$}
\label{fig:ROMinference1c}
\end{subfigure}
\begin{subfigure}[t]{0.47\textwidth}
\includegraphics[width=1\textwidth]{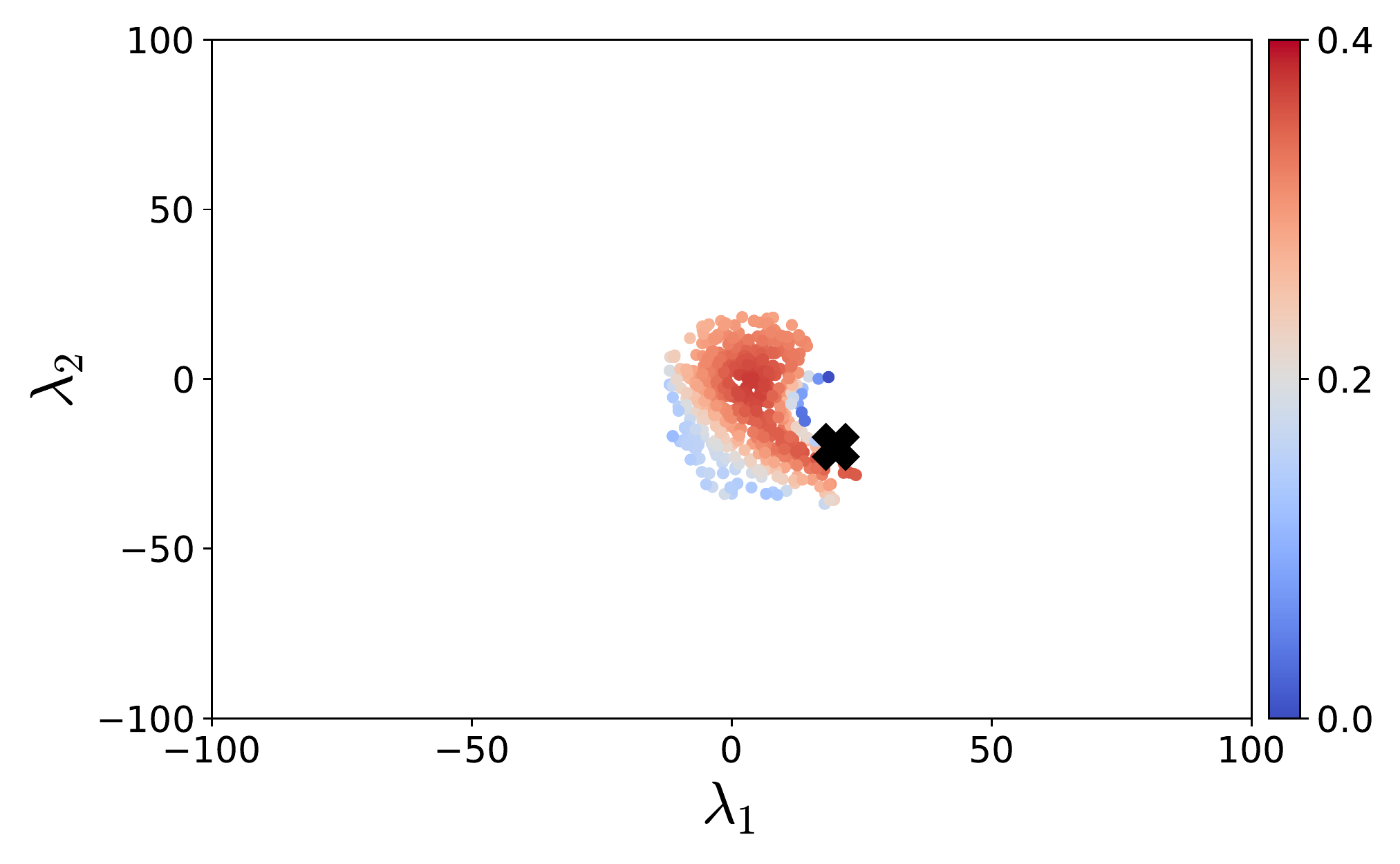}
\caption{7,000 time steps, $\sigma = 2100$}
\label{fig:ROMinference1d}
\end{subfigure}
\caption{Posterior distributions (using a POD-DEIM\new{-Galerkin} ROM as a surrogate for the LDOs) for $\bm{\lambda}$ generated with 1,000 MCMC samples ($\bm{\lambda}^{\text{true}}$ is shown as the black x). Respective values of $\sigma$ used in the likelihood function (Eq.~\ref{eq:Likelihood}) are noted.}
\label{fig:ROMinference1}
\end{figure}

Thus, we are led to conclude that the best compromise available to us
is a middle ground in which we collect just enough observable data
from the ROM \new{in the MCMC algorithm} such that the relevant structural differences begin to
manifest themselves, but not so much that the growing inaccuracies of
the ROM incorrectly skew our inference
results. Fig.~\ref{fig:ROMinference1} provides some empirical evidence
of this. We see that if we limit ourselves to a relatively low number
of snapshots (Fig.~\ref{fig:ROMinference1a}) for inference purposes,
then the highest likelihood samples in our posterior lie along a broad
swath in parameter space. On the other hand, using a relatively high
number of snapshots (Fig.~\ref{fig:ROMinference1d}) collapses the
ridge, but the posterior now completely misses the true
parameters. A number of snapshots between these two extremes
(Fig.~\ref{fig:ROMinference1b} or \ref{fig:ROMinference1c}) produces a
reasonable compromise between accuracy and precision.

\new{We should also briefly discuss a procedure for determining the
  temporal threshold at which we stop collecting data due to the
  decreasing fidelity of the ROM. The most straightforward way to do
  this is to simply choose some random samples in LDO parameter space,
  do a forward simulation of the corresponding full LDO without a ROM,
  and then redo those calculations using the ROM. One could then
  compare the ROM-calculated solutions against the full LDO
  simulations and get an estimate of the length of time over which the
  ROM is able to provide reasonably faithful predictions. Certainly,
  this approach does add some computational overhead, and depending on
  the computational expense of the full LDO system, the user might be
  limited in terms of how many full candidate LDO simulations they can
  afford to do. Notwithstanding, this is a simple and accurate way of
  estimating the ROM’s fidelity over the relevant parameter space. We
  would recommend doing this in practice.}

\new{Additionally, we should add that there are statistical tools for
  explicitly integrating terms that attempt to account for model error
  into computer simulations (see, for
  example,~\cite{kennedy2002bayesian,brynjarsdottir2014learning}). If
  we wanted to go a step further with the ROM-error
  accounting/correction, then it is possible that these tools could be
  useful. One could potentially introduce a discrepancy term that
  would account for the ROM error. Of course, we would have to have
  some prior knowledge about the ROM error in order to do this
  accurately, which we might obtain by the process outlined above.}

Although we have lost some of the statistical precision in the
posterior by using a ROM for the LDO, we have gained a significant
advantage in terms of computational runtime. One way to assess this is
the number of floating point operations required to calculate a LDO
versus a POD-DEIM ROM. For the former, one must first compute the
terms in the basis matrix $\bm{\Psi}$ at each spatial point, which is
a total of $QN_X$ operations (where $N_X$ denotes the number of grid
points). Next, one must compute the matrix product $\bm{\Psi} \bm{P}$,
which requires $\mathcal{O}(N_XQN)$ operations, giving $\mathcal{O}
\left( N_XQN \right)$ operations. In our examples, $N_X = 10^4$, $N =
3$, $M = 5$, $Q = 136$, so $\mathcal{O}(10^6)$ flops are required (for
each timestep). By comparison, the ROM first involves computing the
POD state approximate $\bm{\Phi c}$ at the $d$ DEIM points and at all
of the $M-1$ points in their respective stencils (a total of at most
$Md$ points), which involves $\mathcal{O}(Mdm)$ operations (where $m$
is the number of state POD modes). Next, the LDO must be evaluated on
those $d$ DEIM points, which amounts to calculating the basis elements
and then computing $\bm{\Psi} \bm{P}$ on each of those points, which
involves $ \mathcal{O} \left( QdN \right)$ operations. Finally, the
matrix product $R \bm{f_d}$ must be computed, which involves
$\mathcal{O}(md)$ operations. Altogether, the total number of
operations for each timestep is $\mathcal{O} \left( Mdm + QdN
\right)$. In our examples, we use $m=30$ and $d = 105$, giving a total
of $\mathcal{O}(10^4)$ operations. This is around 100 times fewer
flops and represents a significant computational
alleviation. \new{Empirically, we find that the ROM-based simulation
  is approximately 25 times faster than the full LDO simulation, which
  is reasonably close to our asymptotic operation-count analysis.}

%% \begin{figure}%[h]
%% \centering
%% \includegraphics[width=0.95\textwidth]{ModelVsROMstates3000DT-eps-converted-to.pdf}
%% \caption{Snapshots of the true states compared to ROM-calculated states (after 3,000 timesteps, with $\Delta t = 0.2\frac{1}{\Delta x^2}$).}
%% \label{fig:ActualVsROM3000}
%% \end{figure}

\subsection{Other Issues: Multiple Time Scales, State Space Location of Data}
In this section, we briefly address some questions related to the
effect of the underlying dynamics -- and the numerical code that
implements those dynamics -- on the accuracy of the inferred
LDO. Specifically, we consider how LDO accuracy is affected by
dynamics with fast/slow time scales, and how it is affected by the
location of the data in state space. This is motivated by a dynamical
systems perspective: there could be regions of state space in which
some terms in the governing dynamics are relatively small in magnitude
(e.g., close to a fixed point where linear dynamics dominate), which
might present a problem to the regression method that attempts to
infer them. The practical concern this presents is that in order to
accurately learn those terms, it may be necessary to sample the
dynamics at many timescales and many locations in state space -- which
could be expensive -- in order to infer a reasonable LDO.

We also note that in these examples, we use a set of features that
does not correspond exactly to the operators used in the
discretization of the numerical solver. This is important in the sense
that we would like for our feature set to be both agnostic to and
independent of the discretization used in the actual numerical solver
-- if this were not the case, then we would always have to use a
feature set that included an exact replication of all the algorithms
implemented in the numerical code, which is neither desirable (nor
feasible, often) from the standpoint of system identification.

A nice property of the RSW equations is that they allow for the two
fundamental classes of fluid motion found in oceanic dynamics: slow
vortical motions and fast (inertia-gravity) wave motions. Thus, they
provide a nice testbed for us to investigate the effects of multiple
time scales and the location of the data in state space on LDO
accuracy. To do that, we assess the accuracy of the inferred LDO in
three separate cases. In the first case, the data used for the LDO
inference is a trajectory along the fast mode of the linearized
system. In the second case, the data used is drawn from a trajectory
along the slow linear mode, while the third case is a more general
trajectory that includes both fast and slow components.

\begin{figure}
\centering
\includegraphics[width=0.9\textwidth]{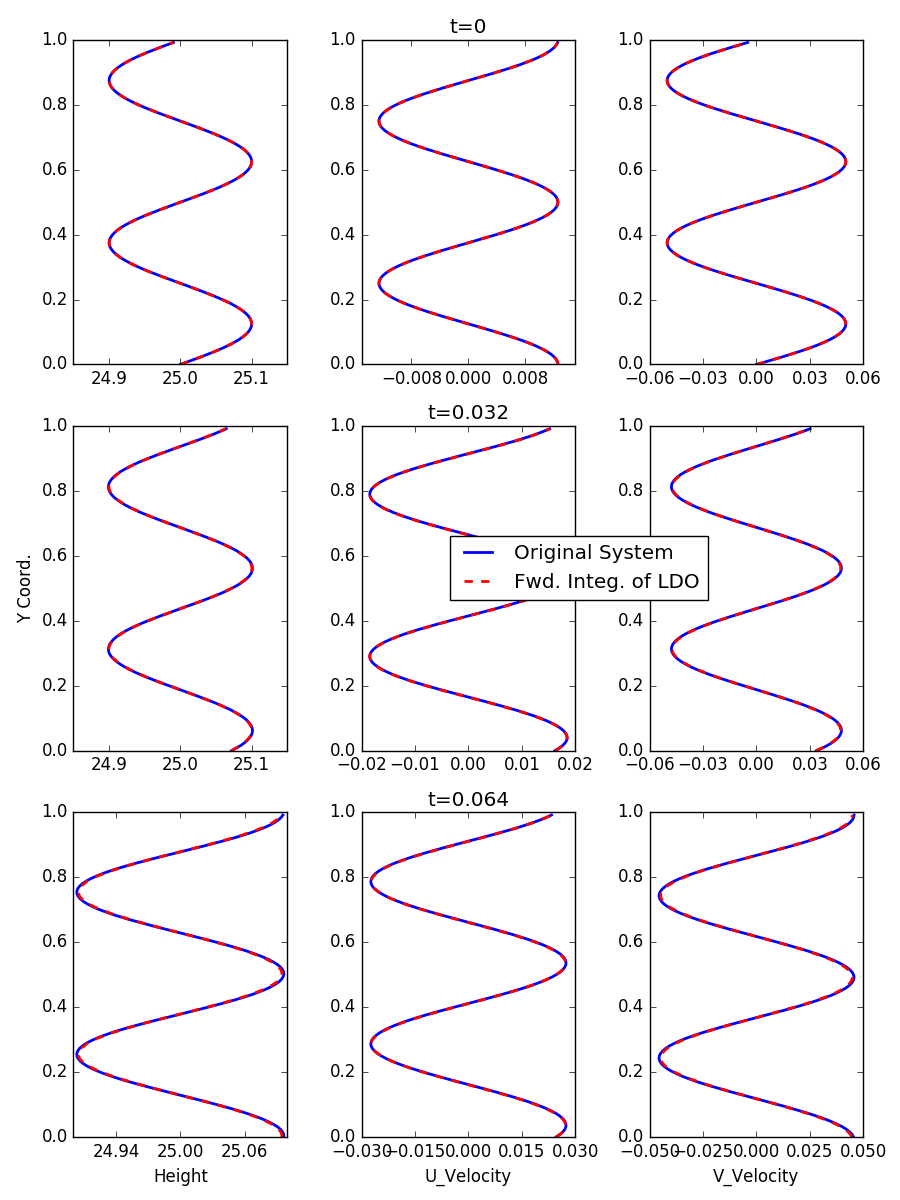}
\caption{Forward integration of the inferred LDO is seen to correctly
  predict (extrapolate) the evolution of fast waves in the system.}
\label{fig:fast}
\end{figure}

Fig.~\ref{fig:fast} shows a comparison of the evolution of the state
of the system as obtained by the forward-integration of the inferred
LDO against the evolution in the original system for a fast mode
initial condition. We note that the LDO inference used data from the
original solver only up to a time of 0.03 so that the snapshots shown
at times 0.032 and 0.064 are both extrapolations from the point of
view of LDO inference. The comparison is seen to be good. Thus, we see
that in order to learn the fast linear dynamics, we need not worry
about collecting data for long periods of time or in multiple regimes
of state space.

As mentioned earlier, we would also like to confirm in these examples
that the LDO inference procedure does not demand that the features
used be exact replicas of those operators implemented in the numerical
solver. Thus, the data used for the LDO inference come from a
pseudo-spectral solver, while the actual features themselves are
consistent with the formulation of~\S[\ref{sec:LDO}]. In light of the
results in Fig.~\ref{fig:fast}, we see that for the fast mode
dynamics, it is indeed possible to infer a reasonable LDO using
features that are independent from the operators implemented in the
numerical solver.

\begin{figure}
\centering
\includegraphics[width=0.9\textwidth]{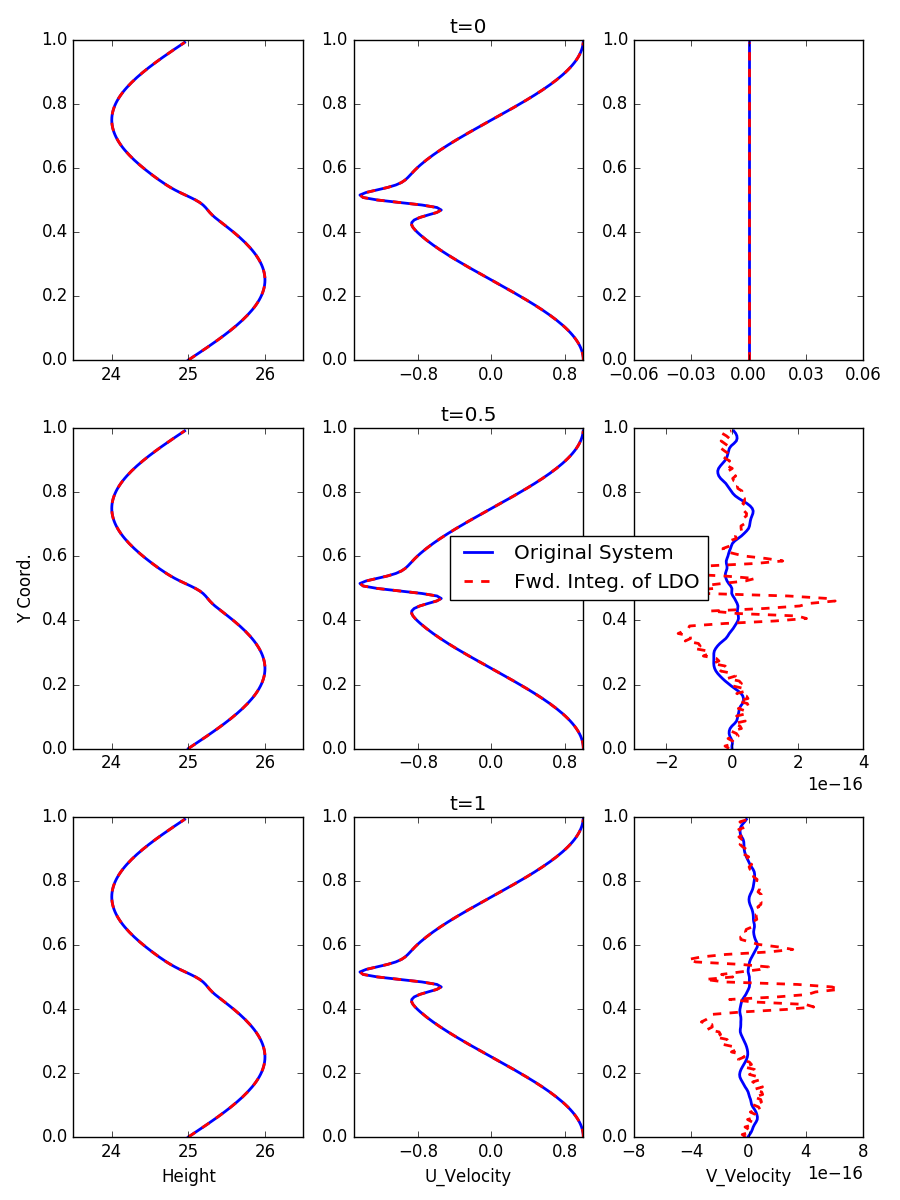}
\caption{As in Fig.~\ref{fig:fast}, the \old{learnt}\new{learned} LDO map is seen to be
  able to capture the slow dynamics correctly as well. Note the
  differences in v-velocity are the level of machine precision
  (round-off errors).}
\label{fig:slow}
\end{figure}

We now consider the evolution of a slow-mode initial
condition. Fig.~\ref{fig:slow} shows a similar comparison as in the
fast-mode initial condition case previously discussed. In this case,
the LDO inference used data from the original solver only up to a time
of 0.01 so that the snapshots shown at times 0.5 and 1.0 are both
(long-term) extrapolations from the point of view of LDO
inference. Again, the comparison is seen to be good, indicating that
-- as in the previous case -- it is not necessary to exhaustively
sample wide swaths of state space for long time periods, or to use a
set of features that correspond exactly to the operators used in the
numerical solver.

\begin{figure}
\centering
\includegraphics[width=0.9\textwidth]{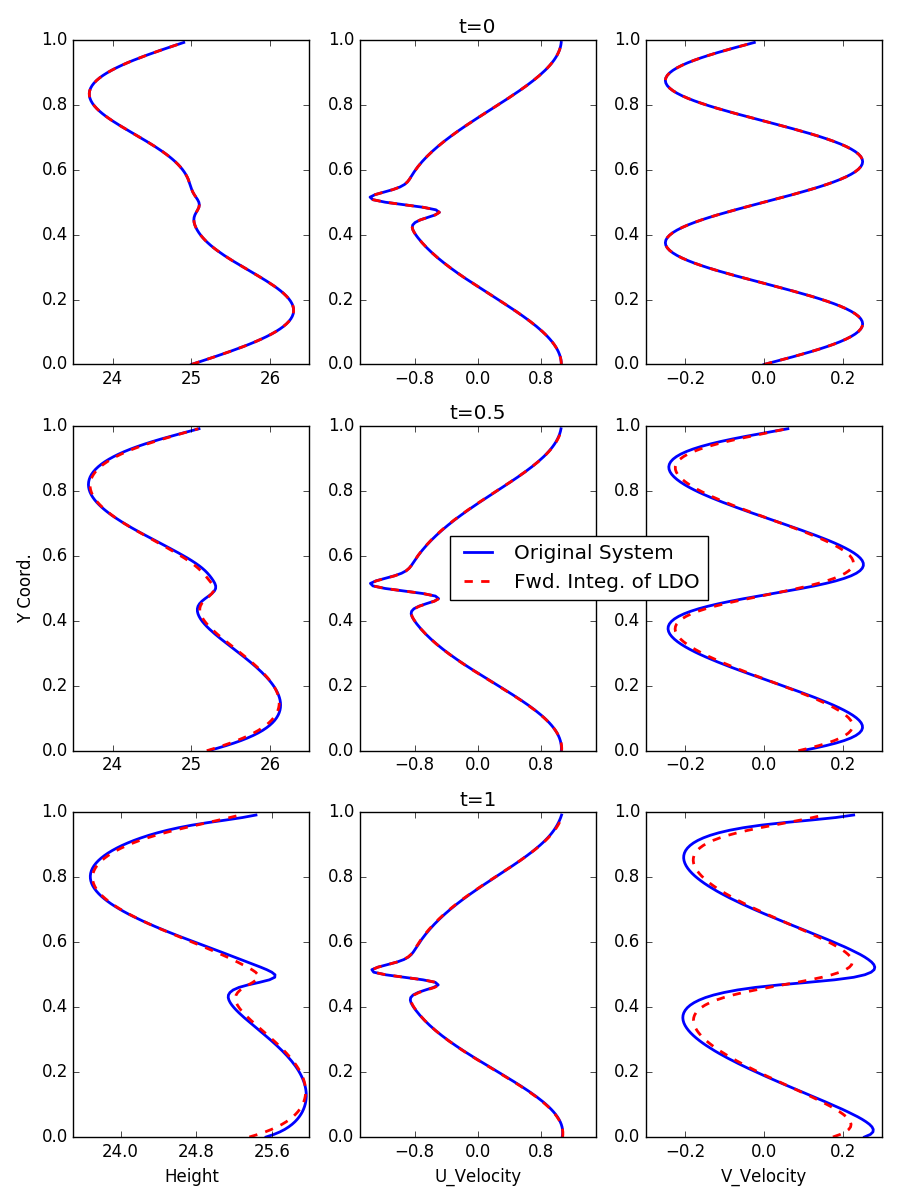}
\caption{In the context of a more generic initial condition that
  involves both fast wave modes and slow vortical modes, significant
  differences are seen in predictions based on the inferred LDO
  suggesting the need for more sophisticated approaches to LDO
  inference.}
\label{fig:mixed}
\end{figure}

Finally, we consider the evolution of an initial condition that
represents a combination of both fast and slow mode components, as
occurs in the oceans and the atmosphere. In this example, the
u-component of velocity and the height fields contain both the slow
vortical components and the fast wave components whereas the
v-component of velocity contains only the fast wave component. This
case represents a qualitatively important difference from the previous
two, since our initial condition no longer lies along an eigenvector
of the linearized system, and so nonlinear terms may be more
pronounced for this case. Fig.~\ref{fig:mixed} displays a comparison
of the evolution of the state of the system as obtained by the
forward-integration of the inferred LDO against the evolution in the
original system. In this example, data up to a time of 0.75 was used
in inferring the LDO, thus making the snapshots at a time of 1.00
LDO-extrapolations. In contrast to the previous cases, the LDO
evolution here is seen to be only approximately correct.

One possible explanation for these discrepancies is that there is some
error in the inferred nonlinear terms of the LDO, which was
effectively irrelevant in the previous two cases where linear dynamics
dominated. This could be the case, for example, if the features we
have used can only at best provide a cruder approximation of the
dynamics relative to the more accurate psuedo-spectral solver
implementation. Another possibility is that the flow evolution is
highly sensitive to location in state space, and so even slight
differences between the two systems compound quickly in time. This
could be the case if, for example, the underlying flow dynamics were
chaotic in this regime of state space. While the former explanation is
an artifact of the inference procedure, the latter is a fundamental
property of the dynamics, and it is not immediately clear which of the
two (or both) accounts for the observed discrepancies. Regardless, in
order to not detract from the main narrative of the article, we
postpone a discussion of such issues and attendant modifications to
the LDO inference procedures to a future article.

\section{Conclusions}
The intended purpose of this paper was to (1) define a mathematical
representation of a model structure that can be manipulated easily and
independently of the original numerical source code, (2) infer this
structure non-intrusively using simulation output, (3) efficiently
approximate the impact of perturbations to a model structure on the
model dynamics, using ROM techniques, and (4) sample the space of
structural uncertainties consistent with data/physical constraints,
and propagate these to uncertainties in predictions (Bayesian
calibration).

The testbed problem on which we demonstrated these objectives was the
2D rotating shallow water equations. We successfully inferred LDO
approximations to the RSW continuum dynamics using two different sets
of features defined on the stencil field values (i.e., quadratic
polynomials and differential operator approximations). We then saw how
it helped to introduce some assumptions about allowable LDO
structures, using energy conservation laws and sparsity-promoting
regression. We constructed a ROM for the purpose of efficient Bayesian
calibration, using non-intrusive system identification
(POD-DEIM-Galerkin projection). We found that the predictive accuracy
POD-DEIM-Galerkin model decreases monotonically with time; this was
expected, since in the Bayesian calibration we are changing the
dynamics from which the ROM was derived. However, we also found some
evidence that the ROM can be sufficiently accurate for an initial
period of time such that we can ``rule out'' large areas of LDO
parameter space before it degrades too much. Using this ROM, we were
able to efficiently generate and search many hypothetical model
structures in LDO parameter space and infer those which were
consistent with a coarse observable quantity of interest.

The envisioned practical setting for these ideas is
quantification/learning of subgrid physics. Our intent is to build a
dictionary of local features that is capable of representing a suite
of subgrid models, learn several such numerical closure models in that
basis, and then search around that area of feature space for a LDO
that reproduces a QOI from some dataset that we regard as the
``truth''.

As we have demonstrated the basics on a toy problem (i.e., the shallow
water equations), the next step is to apply these ideas to more
complex numerical codes. We are currently attempting to do this with
MPAS-Ocean~\cite{MPASOcean}, which is an unstructured code for ocean
dynamics that has several horizontal/vertical mixing parameterizations
and closures implemented. We hope to learn these closure models in a
local feature space, and then use the results to guide an inference of
the subgrid physics for an emprical dataset. This broad topic will be
the subject of future follow-up studies on our part.

\section{Acknowledgments}
Los Alamos Report LA-UR-18-21957. Funded by the Department of Energy
at Los Alamos National Laboratory under contract DE-AC52-06NA25396.

%\clearpage

\bibliographystyle{IJ4UQ_Bibliography_Style.bst}
\bibliography{LDOpaper}

\end{document}